\documentclass{article}
\usepackage[utf8]{inputenc}
\usepackage[margin=2cm]{geometry}
\usepackage{bm} % for bold maths symbols
\usepackage{enumitem} % for more flexibility with lists
\usepackage{array}  % to give vertical centering in tables
\usepackage[pdftex]{graphicx} % for including images
\usepackage{xcolor} % for different text colours
\usepackage{subcaption} % for side by side figures and tables
\usepackage{amssymb} % for mathbb
\usepackage{amsmath} % for align and tfrac
\usepackage{stmaryrd} % for double brackets
\usepackage{authblk} % for author and affiliations
\usepackage[normalem]{ulem} % for strike-outs in editing
\usepackage{fancyhdr} % For headers and footers
\usepackage{lineno} % for line numbers
\usepackage[square,numbers]{natbib}
\newcommand{\grad}[1]{{\bm{\nabla}{#1}}}

\newcommand{\pfrac}[2]{\frac{\partial{#1}}{\partial{#2}}}
\newcommand{\dx}[1]{\hspace{0.75mm}\mathrm{d}{#1}}

\newcommand{\curlybrackets}[1]{\left\lbrace{#1}\right\rbrace}

\fancypagestyle{plain}{
\fancyhf{}
\lfoot{\textcopyright \ Crown Copyright, Met Office}
\cfoot{\thepage}
}

\title{Improving the accuracy of discretisations of the vector transport equation on the lowest-order quadrilateral Raviart-Thomas finite elements}

\author[1]{T. M. Bendall}
\author[2]{G. A. Wimmer}
\affil[1]{Dynamics Research, Met Office, Exeter, UK}
\affil[2]{Los Alamos National Laboratory, Los Alamos, New Mexico 87545, USA}
\date{}

\begin{document}
\maketitle

% ----------------------------------- %
% Abstract

\begin{abstract}
\noindent Within finite element models of fluids, vector-valued fields such as velocity or momentum variables are commonly discretised using the Raviart-Thomas elements.
However, when using the lowest-order quadrilateral Raviart-Thomas elements, standard finite element discretisations of the vector transport equation typically have a low order of spatial accuracy.
This paper describes two schemes that improve the accuracy of transporting such vector-valued fields on two-dimensional curved manifolds. \\
\\
The first scheme that is presented reconstructs the transported field in a higher-order function space, where the transport equation is then solved.
The second scheme applies a mixed finite element formulation to the vector transport equation,  simultaneously solving for the transported field and its vorticity.
An approach to stabilising this mixed vector-vorticity formulation is presented that uses a Streamline Upwind Petrov-Galerkin (SUPG) method.
These schemes are then demonstrated, along with their accuracy properties, through some numerical tests.
Two new test cases are used to assess the transport of vector-valued fields on curved manifolds, solving the vector transport equation in isolation.
The improvement of the schemes is also shown through two standard test cases for rotating shallow-water models. 
\end{abstract}

% Comment / uncomment to add / remove line numbers
% Remember to find+replace to comment / uncomment all % \begin{linenomath} and
% % \end{linenomath} to ensure line numbers don't get disrupted by equations,
% and that text after equations isn't indented
% \linenumbers

% ----------------------------------- %
% Main text

\section{Motivation} \label{sec:motivation}
Many numerical models of fluids involve transporting the velocity or momentum field, via solving an equation such as
% \begin{linenomath}
\begin{equation} \label{eqn:vector transport}
\pfrac{\bm{F}}{t} + \left(\bm{v \cdot \nabla} \right) \bm{F} = \bm{0}.
\end{equation}
% \end{linenomath}
In this work, \eqref{eqn:vector transport} and its variants are referred to as the \textit{vector transport equation}.
In \eqref{eqn:vector transport},  $\bm{v}$ and $\bm{F}$ are vector-valued functions and $\bm{F}$ is transported by $\bm{v}$.
A major class of such numerical models are those that use finite element methods, which have a major advantage of being easy to formulate on arbitrary meshes.
In finite element methods, fields are expressed as the sum of a finite number of basis functions multiplied by coefficients.
Each basis function is localised to a cell or a small number of cells on the mesh.
The choice of basis functions and their continuity between cells is typically referred to as the \textit{finite element}.
\\
\\
As argued by \cite{bernard2009high}, finite element discretisations can also offer advantages when considering the transport of vectors on two-dimensional curved manifolds, such as the surface of the sphere.
In such cases, the vector transport equation generally includes \textit{metric terms}, which describe accelerations induced by the curvature of the manifold itself.
As discussed by \cite{bernard2009high}, there are two standard approaches to handling these terms.
Firstly the metric terms can be explicitly included in the equation, which typically involves evaluating the Christoffel symbols describing the curvature of the manifold.
However for a general manifold this evaluation may not be straightforward, which can make this approach difficult or even impossible.
In the alternative approach, the vectors are described in three Cartesian components, which adds an extra dimension to the equation.
Then no metric terms appear in the transport equations for the components, but a third unknown has been added, and a constraint must also be applied to keep the transported vector in the tangent bundle of the manifold.
Instead, finite element discretisations can combine the benefits of these two approaches.
By writing the equation in a weak integral form and numerically evaluating the integral in Cartesian coordinates, the explicit evaluation of Christoffel symbols can be avoided.
At the same time, no third component is added and the transported vector will naturally be tangent to the manifold.
\\
\\
One family of finite elements that is often used for velocity or momentum variables is the Raviart-Thomas family, which can be defined on triangular or quadrilateral cells.
The previous decade has seen particular interest in these finite elements from the numerical weather prediction (NWP) community.
In many of the finite difference or finite volume methods used historically by this community, the density/pressure and velocity variables have been staggered according to the Arakawa C-grid of \cite{winninghoff1968,mesinger1976numerical,arakawa1977computational},
due to its good representation of the wave modes of the shallow-water equations \cite{staniforth2012horizontal}.
It was shown by \cite{cotter2012mixed} that certain choices of finite element pairs for the density/pressure and velocity variables can still replicate the desirable dispersion properties of the Arakawa C-grid in a mixed finite element model of the shallow-water equations.
In particular, the finite element equivalent of the velocity staggering used in the Arakawa C-grid on quadrilateral cells is the lowest-order Raviart-Thomas elements (i.e.\ the elements using the lowest degree polynomials in the basis functions).
Maintaining this equivalence by using the lowest-order finite element spaces can be advantageous for other reasons, for instance that it can simplify the coupling to parametrisations which are used to represent unresolved physical processes in NWP models. \\
\\
It is these properties that have seen the lowest-order Raviart-Thomas elements become candidates for use in NWP models.
For instance, the UK Met Office will use them in its next-generation model, LFRic\footnote{For more information on LFRic, see \cite{adams2019lfric} and \cite{melvin2019mixed}.} (named after Lewis Fry Richardson).
The Met Office currently uses a longitude-latitude grid for its global simulations, which suffers from the \textit{pole problem}.
The convergence of meridians at the poles of the grid has begun to lead to bottlenecks in data communication on massively parallel supercomputers.
The scalability of the model is then compromised, and the current forecasting model will be unable to exploit the computational power of the next-generation of supercomputers.
Moving to a finite element formulation facilitates the move to a cubed-sphere grid, which is quasi-uniform over the sphere and should avoid these scalability bottlenecks.\\
\\
However, a challenge with using these lowest-order elements is that typical discretisations of the transport equation do not have a satisfactory order of accuracy with respect to the grid spacing, (as discussed by \cite{staniforth2012horizontal} this should be at least approaching second-order). 
One route to circumvent this is to use higher-order finite difference or finite volume methods to build up transport stencils, which is the approach used by \cite{melvin2019mixed}.
Unfortunately finite difference or finite volume methods are not supported in many finite element software systems, where these methods may then be unfeasible.
In any case, as argued by \cite{bernard2009high}, finite element discretisations may offer particular advantages for discretising the vector transport equation.
The motivation is then to find finite element discretisations that do deliver improved accuracy, while in this work the computational cost of such schemes are of secondary concern. 
\\
\\
One approach to tackle the low order of accuracy of transport schemes for the lowest-order elements was presented by \cite{bendall2019recovered}, which introduced a method of recovering fields in a higher-order finite element space for solving the transport equation.
This resulted in higher-order accuracy overall while using the lowest-order finite elements.
However \cite{bendall2019recovered} did not present a method that could be used for transporting vector-valued fields on curved manifolds.
In this paper, two discretisations of the vector transport equation are shown to improve the order of accuracy for transport with the lowest-order Raviart-Thomas elements on quadrilateral cells, and crucially when the transport is on a two-dimensional curved manifold.
The first method extends the recovery approach of \cite{bendall2019recovered} to curved manifolds, reconstructing the vector-valued field in a higher-order function space.
The second method adapts a mixed finite element formulation similar to that of \cite{bauer2018energy} to the vector transport equation.
This scheme simultaneously solves for the transported vector and its vorticity.
A stabilisation based on a Streamline Upwind Petrov-Galerkin (SUPG) approach is then presented for this scheme, based on \cite{wimmer2020energy}.
\\
\\
The remainder of the paper is laid out as follows.
In Section \ref{sec:background}, some background is given to the vector transport equation and the Raviart-Thomas elements, alongside a standard upwind finite element scheme for \eqref{eqn:vector transport} that is used as a benchmark.
The two schemes that improve on this are described in Sections \ref{sec:recovery} and \ref{sec:vorticity} respectively.
Section \ref{sec:recovery} reviews the recovered transport approach of \cite{bendall2019recovered} before extending it for the Raviart-Thomas elements,
while Section \ref{sec:vorticity} describes a mixed vector-vorticity discretisation like that of \cite{bauer2018energy} in the context of the vector transport equation, and presents the new SUPG stabilisation to it.
The new schemes are demonstrated through some test cases in Section \ref{sec:results}, that cover both the vector transport equation on its own and also within a shallow-water model.

\section{Background} \label{sec:background}

\subsection{The Vector Transport Equation} \label{sec:forms of eqn}
This work considers the transport of some vector-valued field $\bm{F}(\bm{x},t)$ by some other vector-valued field $\bm{v}(\bm{x},t)$, where $\bm{x}$ is the position vector in the domain $\varOmega$ and $t$ is the point in time.
The domain $\varOmega$ is a two-dimensional differentiable manifold, which may be embedded in two-dimensional space (so that the domain is a plane) or in three-dimensional space (for instance when the domain is the surface of a sphere).
The vectors $\bm{v}$ and $\bm{F}$ live in the tangent bundle of $\varOmega$, and so can be expressed locally through two scalar components. \\
\\
This section briefly considers different forms of the vector transport equation.
It is most easily expressed as \eqref{eqn:vector transport}, which is referred to as the \textit{advective form}, and which is repeated again here:
% \begin{linenomath}
\begin{equation} \label{eqn:vec transport -- advective}
\pfrac{\bm{F}}{t} + \left(\bm{v \cdot \nabla} \right) \bm{F} = \bm{0}.
\end{equation}
% \end{linenomath}
Before writing the next form of the equation, we introduce the perpendicular operation, denoted by superscript $^\perp$.
Defining the unit normal outward from the manifold as $\widehat{\bm{N}}$, then $\bm{F}^\perp$ is given by
% \begin{linenomath}
\begin{equation} \label{def:perp}
\bm{F}^\perp := \widehat{\bm{N}}\times\bm{F},
\end{equation}
% \end{linenomath}
with $\times$ denoting the cross product.
If $\varOmega$ is a plane with components labelled $x$ and $y$, then  $\bm{F}^\perp=\left(-F_y,F_x\right)$.
With this definition, an alternative form of the vector transport equation is
% \begin{linenomath}
\begin{equation}\label{eqn:vec transport -- vec inv}
\pfrac{\bm{F}}{t} + \left(\bm{\nabla}^\perp\bm{\cdot F}\right) \bm{v}^\perp
+ \frac{1}{2}\grad{(\bm{v\cdot F})}
+ \frac{1}{2}\left[\left(\bm{\nabla} \bm{F}\right)\bm{\cdot v}     
-\left(\bm{\nabla}\bm{v}\right)\bm{\cdot F}  \right] = \bm{0},
\end{equation}
% \end{linenomath}
where $\bm{\nabla}$ applied to a vector is tensor-valued, and we call the terms featuring this the \textit{vector-gradient} terms.
The second term of \eqref{eqn:vec transport -- vec inv} might be more easily recognised as the two-dimensional version of $\left(\bm{\nabla}\times\bm{F}\right)\times\bm{v}$. \\
\\
If the \textit{vorticity} is defined by
% \begin{linenomath}
\begin{equation} \label{def:vorticity}
\zeta := \bm{\nabla}^\perp\bm{\cdot F},
\end{equation}
% \end{linenomath}
then \eqref{eqn:vec transport -- vec inv} can also be written in \textit{vorticity form}:
% \begin{linenomath}
\begin{equation}\label{eqn:vec transport -- vorticity}
\pfrac{\bm{F}}{t} + \zeta \bm{v}^\perp
+ \frac{1}{2}\grad{(\bm{v\cdot F})}
+ \frac{1}{2}\left[\left(\bm{\nabla} \bm{F}\right)\bm{\cdot v}     
-\left(\bm{\nabla} \bm{v}\right)\bm{\cdot F}  \right] = \bm{0}.
\end{equation}
% \end{linenomath}
The velocity $\bm{u}$ in fluid dynamics models is self-transporting, and in that context the vector transport equation becomes a form of the Burgers' equation.
With $\bm{v}=\bm{F}=\bm{u}$, the vector-gradient terms of \eqref{eqn:vec transport -- vorticity} cancel to yield the \textit{vector-invariant} form:
% \begin{linenomath}
\begin{equation} \label{eqn:vector invariant}
\pfrac{\bm{u}}{t} + \zeta \bm{u}^\perp + \tfrac{1}{2}\bm{\nabla}\left(\bm{u}\bm{\cdot}\bm{u}\right)= \bm{0}.
\end{equation}
% \end{linenomath}
Finally, although not considered in this work, the vector transport equation can also be written in \textit{flux form}:
% \begin{linenomath}
\begin{equation} \label{eqn:vec transport -- flux}
\pfrac{\bm{F}}{t} + \bm{\nabla\cdot}\left(\bm{v}\otimes\bm{F}\right)
- \left(\bm{\nabla\cdot v}\right)\bm{F}=\bm{0},
\end{equation}
% \end{linenomath}
where $\otimes$ denotes the outer product of two vectors.
\subsection{Raviart-Thomas Elements}
The Raviart-Thomas family is an important class of finite elements used to describe two-dimensional vector fields.
These elements were introduced for both triangular and quadrilateral elements by \cite{raviart1977mixed}, who used them to solve the Poisson equation with a mixed finite element discretisation.
The Raviart-Thomas elements come in two varieties: those that preserve the normal components of vectors between cells, and those that preserve the tangential components.
These former elements (preserving normal components) are known as $H(\mathrm{div})$-conforming as functions in these elements have square-integrable divergence. 
The latter (preserving tangential components) are $H(\mathrm{curl})$-conforming, which in the context of a two-dimensional manifold means that for all fields $\bm{u}$ in the corresponding finite element space, $\bm{\nabla}^\perp\bm{\cdot u}$ is square-integrable (as well as $\bm{u}$ also being square-integrable).
As discussed by \cite{arnold2015finite}, the quadrilateral Raviart-Thomas elements can be represented as the tensor-product of one-dimensional elements.
For definitions and more thorough descriptions of the elements, see \cite{brezzi2012mixed}.\\
\\
This work uses the nomenclature of \cite{arnold2014periodic},
with $\mathrm{RTc}^e_k$ representing the $k$-th order $H(\mathrm{curl})$-conforming elements on quadrilateral elements and $\mathrm{RTc}^f_k$ representing the $H(\mathrm{div})$-conforming elements on quadrilateral elements.\\
\\
The field of finite element exterior calculus (see \cite{arnold2010finite}) explains how some finite elements can be related to others through the action of the exterior derivative.
Such spaces can be part of a de Rham complex, which is the chain of spaces obtained by application of the exterior derivative.
For instance, taking the divergence of a field in the $\mathrm{RTc}^f_k$ space yields a field in the discontinuous Galerkin space $\mathrm{DG}_{k-1}$.
The main families of finite elements that form de Rham complexes are captured in the periodic table of finite elements \cite{arnold2014periodic}.\\
\\
In a \textit{compatible} finite element model, variables in the discretisation are chosen to lie in the spaces of the discrete de Rham complex that correspond to their continuous analogues.
In this structure, the discrete differential operators preserve vector calculus identities such as $\bm{\nabla}\times\bm{\nabla}f=\bm{0}$ for all scalar $f$.
Applied to fluid dynamics, this structure suggests that the velocity should lie in a $H(\mathrm{div})$-conforming space such as the $\mathrm{RTc}^f_k$ elements.
As mentioned in Section \ref{sec:motivation}, it was shown by \cite{cotter2012mixed} that the choice of $\mathrm{RTc}^f_k$-$\mathrm{DG}_{k-1}$ for the wind and height fields on quadrilateral elements in a shallow-water model gives a discretisation with many desirable properties.
As explained by \cite{cotter2012mixed}, this pair of elements has the optimal ratio of degrees of freedom (DoFs) for capturing the wave modes of the shallow-water equations, and it mimics the properties of the popular C-grid staggering used in finite difference models.
For this reason, a similar discretisation will be used in the Met Office's new LFRic model, with the wind lying in the three-dimensional form of $\mathrm{RTc}^f_1$ space.
Although the use of the corresponding Raviart-Thomas elements on triangular elements in a shallow-water model has been investigated elsewhere, notably by \cite{rostand2008raviart} as part of the $\mathrm{RT}^f_k$-$\mathrm{DG}_{k-1}$ pair, as shown by \cite{cotter2012mixed} it suffers from inferior representation of the shallow-water wave modes.
Given this result, this work focuses on quadrilateral elements.

\subsection{An Upwind Finite Element Discretisation} \label{sec:benchmark scheme}
This section describes a simple finite element discretisation for the advective form \eqref{eqn:vec transport -- advective}, which is used in the results of Section \ref{sec:results} as a benchmark to compare with the improved schemes of Sections \ref{sec:recovery} and \ref{sec:vorticity}.
This discretisation is a generalisation of the upwind discontinuous-Galerkin method (first used by \cite{reed1973triangular}) to vector-valued fields. \\
For an overview of these methods see for instance \cite{cockburn2001runge}.
For these methods, the time discretisation generally does not affect the spatial accuracy, so discussion of the time discretisation is left until Section \ref{sec:results}.
\\
\\
Let $\bm{v}$ and $\bm{F}$ lie in function space $V_F$ made up of Raviart-Thomas elements.
Multiplying \eqref{eqn:vec transport -- advective} by a test function $\bm{\gamma}\in V_F$, integrating over the domain $\varOmega$ and then integrating by parts gives, $\forall \bm{\gamma}\in V_F$,
% \begin{linenomath}
\begin{equation} \label{eqn:dg upwind}
\begin{split}
\int_\varOmega \bm{\gamma \cdot} \pfrac{\bm{F}}{t} \dx{x}
+ \int_\varGamma \left(\bm{v}^+ \bm{\cdot}\widehat{\bm{n}}^+\right)
\left\llbracket\bm{\gamma}\right\rrbracket_{+}\bm{\cdot}\bm{F}^\dagger \dx{S}
- \int_\varOmega \bm{F\cdot}\left[\bm{\nabla\cdot}\left(\bm{\gamma}\otimes\bm{v}\right)\right]\dx{x} \\
+ \int_\varGamma \left(\bm{v}^+ \bm{\cdot}\widehat{\bm{n}}^+\right)
\left(\bm{F}^\dagger\bm{\cdot}\widehat{\bm{n}}^\dagger \right)
\left(\bm{\gamma}^\ddagger\bm{\cdot}\left[\widehat{\bm{n}}^++\widehat{\bm{n}}^- \right]\right) \dx{S} = 0.
\end{split}
\end{equation}
% \end{linenomath}
Here $\varGamma$ is the set of all interior facets of the domain.
Each side of these facets can be arbitrarily labelled with a $+$ or $-$, and $\widehat{\bm{n}}^+$ is defined as the outward normal from the $+$ side of a facet.
The double square brackets $\left\llbracket\cdot\right\rrbracket_+$ denote the jump of some field over a facet, so that
% \begin{linenomath}
\begin{equation} \label{def:jump}
\left\llbracket \bm{\gamma} \right\rrbracket_+ := \bm{\gamma}^+ - \bm{\gamma}^-.
\end{equation}
% \end{linenomath}
The upwind value at a facet is denoted by the dagger $^\dagger$, and is given by
% \begin{linenomath}
\begin{equation} \label{def:upwind}
\bm{F}^\dagger := \left\lbrace 
\begin{matrix}
\bm{F}^+ & \mathrm{if} \ \bm{v^+ \cdot}\widehat{\bm{n}}^+ \geq 0, \\
\bm{F}^- & \mathrm{if} \ \bm{v^+ \cdot}\widehat{\bm{n}}^- < 0.
\end{matrix}\right.
\end{equation}
% \end{linenomath}
The final term of \eqref{eqn:dg upwind} is a correction to project the upwind term into the tangent bundle,
with the double dagger $^\ddagger$ denoting the downwind term (i.e.\ from the opposite side of the facet to the upwind term).
This correction is similar to that used by \cite{bernard2009high}, so that both sides of the $\left\llbracket \bm{\gamma}\right\rrbracket_+$ term are evaluated in the tangent space of $\bm{F}^\dagger$, on the upwind side of the facet.
In a Cartesian plane, $\widehat{\bm{n}}^+=-\widehat{\bm{n}}^-$ and the correction vanishes, but this is not generally true for a curved manifold.
\\
\\
Although this benchmark scheme performs well for general Raviart-Thomas spaces, it has low-order accuracy for the lowest-order spaces, which is demonstrated in Section \ref{sec:results}.
This poor performance can be understood heuristically by considering the components and basis functions of $\bm{F}$.
Those that are parallel to $\bm{v}$ are linear in a cell in the direction of $\bm{v}$. 
However the components of $\bm{F}$ that are perpendicular to $\bm{v}$ are only constant in a cell in the direction of $\bm{v}$.
Conventional finite element discretisations of spatial derivatives for piecewise constant fields have only first-order accuracy or worse.
For this reason, similar upwind discretisations (such as that of \cite{natale2018variational}) will also suffer from low orders of accuracy when applied to alternative forms of the transport equation such as \eqref{eqn:vec transport -- vec inv} or \eqref{eqn:vec transport -- flux}.

\subsection{Discretisation of the Shallow-Water Equations}\label{sec:shallow-water background}
Some discretisations of geophysical fluids include a step in which the vector transport equation is solved in isolation.
One such discretisation is the shallow-water model of \cite{shipton2018higher} and \cite{gibson2019compatible}, which uses a compatible finite element framework.
This section briefly describes this model, which is used for the demonstrations in Section \ref{sec:results}, applied to the lowest-order finite element spaces that correspond to those used by the Met Office's LFRic model \cite{melvin2019mixed}, so that the velocity field $\bm{u}$ and the depth field $h$ are in the $\mathrm{RTc}^f_1$ and $\mathrm{DG}_0$ spaces respectively.
\\
\\
The rotating shallow-water equations can be expressed as
% \begin{linenomath}
\begin{subequations} \label{eqn:shallow-water}
\begin{align}
& \pfrac{\bm{u}}{t}+(\bm{u\cdot \nabla})\bm{u}+f\bm{u}^\perp+g\grad{(h+h_b)}=0, \label{eqn:shallow-water-momentum}\\
& \pfrac{h}{t} + \bm{\nabla\cdot}(h\bm{u})=0,
\end{align}
\end{subequations}
% \end{linenomath}
where $h_b$ is the height of the lower surface, $f$ is the Coriolis parameter and $g$ is acceleration due to gravity.
For more discussion of the shallow-water equations, see for instance \cite{williamson1992standard}.
The shallow-water model of \cite{shipton2018higher} and \cite{gibson2019compatible} discretises \eqref{eqn:shallow-water} with a time stepping structure that follows the semi-implicit scheme used by both the Met Office's current ENDGame \cite{walters2017met} and new GungHo dynamical cores \cite{melvin2019mixed}.
In this semi-implicit scheme, a time step consists of an outer loop in which the transport terms are evaluated, and an inner loop in which the implicit terms are obtained by solving a linearised form of \eqref{eqn:shallow-water}.
This linear problem is iterated to obtain the variables at the next time step, $\bm{u}^{n+1}$, $h^{n+1}$.
For the linear solver, the hybridised finite element technique presented by \cite{gibson2019compatible} is used.
A thorough description of the semi-implicit time stepping scheme is also given by \cite{gibson2019compatible}.\\
\\
As part of the outer loop of the time step, the transport terms are evaluated from
% \begin{linenomath}
\begin{subequations}
\begin{align}
& \pfrac{\bm{u}}{t} = - (\bm{u}_a\bm{\cdot\nabla})\bm{u}^*,  \label{eqn:sw_vector_transport} \\
& \pfrac{h}{t} = - \bm{\nabla\cdot}(h^n\bm{u}_a),
\end{align}
\end{subequations}
% \end{linenomath}
where $\bm{u}^\ast$ is $\bm{u}^n$ incremented by the explicit pressure gradient and Coriolis terms.
The transporting velocity is $\bm{u}_a=\frac{1}{2}\left(\bm{u}^n+\bm{u}^{(k)}\right)$, where $\bm{u}^{(k)}$ is the latest approximation of $\bm{u}^{n+1}$.
So comparing \eqref{eqn:sw_vector_transport} with \eqref{eqn:vec transport -- advective}, $\bm{u}_a$ plays the role of $\bm{v}$ and $\bm{u}^*$ plays the role of $\bm{F}$.
In this context, the discretisation of \eqref{eqn:sw_vector_transport} can be treated as a ``black box'', in which different schemes for solving the vector transport equation can be used.

\section{Recovery} \label{sec:recovery}
Inspired by the recovered finite element methods of \cite{georgoulis2018recovered} and the embedded transport scheme of \cite{cotter2016embedded}, \cite{bendall2019recovered} presented a transport scheme to improve the spatial accuracy of discretisations for the lowest-order finite element spaces.
This was particularly motivated for the discontinuous Galerkin space of piecewise constants, $\mathrm{DG}_0$.
The recovered transport scheme of \cite{bendall2019recovered} involves recovering the field to be transported from $\mathrm{DG}_0$ to $\mathrm{DG}_1$, and solving the transport equation in $\mathrm{DG}_1$ before projecting the solution back to $\mathrm{DG}_0$.
The recovery is performed using a simple averaging operator, which was indicated by \cite{georgoulis2018recovered} to have second-order accuracy.
If then the transport scheme used for the $\mathrm{DG}_1$ field has second-order accuracy, the whole scheme has second-order accuracy.
However, \cite{bendall2019recovered} focused on Cartesian domains where the velocity field was transported by decomposing the field into orthogonal components and transporting each of these separately.
This recovered approach was also used for solving the transport equation by \cite{bendall2020compatible} in the context of a moist compressible Euler model, but this also only focused on Cartesian domains.
In this section, after reviewing the approach of \cite{bendall2019recovered}, we present the extension to this to achieve higher-order transport on curved manifolds.
\subsection{Review of recovered transport} \label{sec:recovery review}
To start, we define a series of function spaces $\curlybrackets{V_L, \widehat{V}_L, V_R, V_H}$ and operators $\curlybrackets{\mathcal{I}, \mathcal{R}, \mathcal{P}_L, \widehat{\mathcal{P}}_L, \mathcal{J}, \mathcal{T}}$ that are used in the recovered transport scheme of \cite{bendall2019recovered}.
These are summarised in Table \ref{tab:rec_spaces_ops}.
Here a different terminology is used to that of \cite{bendall2019recovered} to make the new scheme clearer. \\
\\
The lowest-order function space is given by $V_L$.
This is the native function space of the transported variable $q$,  so that $q\in V_L$.
The \textit{broken} (fully-discontinuous) form of $V_L$ is then denoted by $\widehat{V}_L$.
The higher-order space in which the transport will happen is $V_H$, and an intermediate space into which $q$ is recovered is $V_R$.
In \cite{bendall2019recovered}, the spaces were chosen so that $V_L\subset V_H$, $\widehat{V}_L\subset V_H$ and $V_R\subset V_H$,
while $V_R$ was assumed to be fully-continuous.
Section \ref{sec:recovery new scheme} relaxes some of these requirements. \\
\\
As in \cite{cotter2016embedded}, an \textit{injection operator} $\mathcal{I}_H:V\to V_H$ identifies a field in one of $V_L$, $\widehat{V}_L$ and $V_R$ as also being a member of $V_H$.
The operator $\mathcal{P}_L:V\to V_L$ is a Galerkin projection from some space into the lower-order space.
Taking arbitrary $\bm{u}\in V$ (e.g.\ $V_H$), the action of $\mathcal{P}_L$ so that $\bm{y} = \mathcal{P}_L\bm{u}$ with $\bm{y}\in V_L$ is given by
% \begin{linenomath}
\begin{equation}
\int_\varOmega \bm{\gamma \cdot y} \dx{x} = \int_\varOmega \bm{\gamma \cdot u} \dx{x}, \ \ \ \ \forall \bm{\gamma} \in V_L.
\end{equation}
% \end{linenomath}
Similarly, $\widehat{P}_L:V_R\to \widehat{V}_L$ is a Galerkin projection.
The key operator in the reconstruction is $\mathcal{R}:V_L \to V_R$.
This is the \textit{recovery operator}, and should have second-order spatial accuracy.
This can be achieved by using an averaging operator, so that for the DoFs of $V_R$ that are shared between cells, the field values are the average of values from neighbouring cells of the field in $V_L$.
At any domain boundaries, improved accuracy can be obtained by extrapolating from values on the interior
(for more details see \cite{bendall2019coupling}). \\
\\
The full operator for reconstructing the higher-order field is $\mathcal{J}:V_L\to V_H$, defined by
% \begin{linenomath}
\begin{equation}
\mathcal{J} := \mathcal{I}_H+ \mathcal{I}_H\mathcal{R} - \mathcal{I}_H\widehat{\mathcal{P}}_L\mathcal{R}.
\end{equation}
% \end{linenomath}
The addition of $\mathcal{I}_H-\mathcal{I}_H\widehat{\mathcal{P}}_L\mathcal{R}$ ensures that the mass of $\bm{u}\in V_L$ and $\mathcal{J}\bm{u}\in V_H$ will be preserved in each cell, and that $\mathcal{P}_L\mathcal{J}\bm{u}=\bm{u}$, so that if no transport happens then $\bm{u}$ will remain unchanged.\\
\\
Finally, $\mathcal{T}:V_H\to V_H$ performs a single transport step in $V_H$.
Denoting the value of $\bm{u}\in V_L$ at the $n$-th time step as $\bm{u}^n$, a whole transport step is described by
% \begin{linenomath}
\begin{equation} \label{eqn:recovered step}
\bm{u}^{n+1} = \mathcal{P}_L\mathcal{T}\mathcal{J}\bm{u}^{n}.
\end{equation}
% \end{linenomath}

\begin{table}
\centering
\begin{subtable}[t]{.45\linewidth}
\centering
\caption{Nomenclature of function spaces}
\begin{tabular}{c|l}
\textbf{Space} & \textbf{Description} \\
\hline
$V_L$ & Native lower-order space of transported variable     \\
$V_H$ & Higher-order space to perform transport in           \\
$\widehat{V}_L$ & Broken (fully-discontinuous) form of $V_L$ \\
$V_R$ & Higher-order space to recover into \\
$\widehat{V}_R$ & Broken (fully-discontinuous) form of $V_R$
\end{tabular}
\label{tab:space defs} 
\end{subtable}
\hfill
\begin{subtable}[t]{.45\linewidth}
\centering
\caption{Nomenclature of operators}
\begin{tabular}{l|l}
\textbf{Operator} & \textbf{Description} \\
\hline
$\mathcal{I}_H:V \to V_H$  & Injection operator into $V_H$ \\
$\mathcal{P}_L:V \to V_L$ & Projection operator into $V_L$ \\
$\widehat{\mathcal{P}}_L:V_R \to \widehat{V}_L$ & Projection operator into $\widehat{V}_L$ \\
$\mathcal{P}_H:V\to V_H$ & Projection operator into $V_H$ \\
$\widehat{\mathcal{P}}_R:V_L\to \widehat{V}_R$ & Projection operator into $\widehat{V}_R$ \\
$\mathcal{A}:\widehat{V}_R\to V_R$ & Averaging operator \\
$\mathcal{R}:V_L \to V_R$ & Recovery operator  \\
$\mathcal{J}:V_L \to V_H$  & Full reconstruction operator  \\
$\mathcal{T}:V_H \to V_H$ & Transport operator in $V_H$
\end{tabular}
\label{tab:operators}
\end{subtable}
\caption{A summary of the variables used to describe the function spaces and operators involved in the recovered transport scheme. The space $V$ without a subscript is used to represent a range of the other defined spaces.}
\label{tab:rec_spaces_ops}
\end{table}

\subsection{Extension to curved manifolds} \label{sec:recovery new scheme}
In \cite{bendall2019recovered}, the higher-order space $V_H$ for the transport of scalar-valued fields was taken as the $\mathrm{DG}_1$ space, while $V_R$ was the linear continuous Galerkin space $\mathrm{CG}_1$.
To transport the velocity field, it was separated into orthogonal components which were each separately reconstructed in $\mathrm{DG}_1$.
This was possible because \cite{bendall2019recovered} only considered Cartesian domains.\\
\\
However using this approach on curved manifolds presents problems.
For instance, consider two vectors at two different points of the manifold, both pointing along the geodesic that joins the two points.
These vectors will generally lie in two different tangent planes.
A vector lying on the midpoint of the geodesic between the two points should not be reconstructed by using the average of the Cartesian components.
This would generally not lie in the tangent space itself, and its projection into the tangent space will likely under-approximate the vector's size.
For some domains this could be resolved by averaging in some other orthogonal coordinate system, but this work is motivated by geophysical applications and in particular the sphere, where the topology also presents challenges (for instance when a spherical-polar coordinate system is used then the components do not make sense at the poles).
In this section, we extend the scheme of \cite{bendall2019recovered} by careful choice of the spaces and operators described in Section \ref{sec:recovery review} so as to avoid these problems and achieve a higher-order transport scheme for velocities in the lowest-order Raviart-Thomas spaces.\\
\\
The broad structure of the scheme is the same, following equation \eqref{eqn:recovered step}, so that
% \begin{linenomath}
\begin{equation}
\bm{u}^{n+1} = \mathcal{P}_L\mathcal{T}\mathcal{J}\bm{u}^{n}.
\end{equation}
% \end{linenomath}
The main difference is that the operator $\mathcal{J}$ will be defined differently.
Although it is still required that that $V_L\subset V_H$, it is no longer assumed that $V_R$ or $\widehat{V}_L$ are subsets of $V_H$.
Another difference is the introduction of $\widehat{V}_R$, the space of broken elements of $V_R$. \\
\\
The recovery operator $\mathcal{R}:V_L\to V_R$ is split into two steps: firstly a Galerkin projection $\widehat{P}_R:V_L\to \widehat{V}_R$, and secondly an averaging operator $\mathcal{A}:\widehat{V}_R\to V_R$.
The averaging operator restores the continuity of a field in $\widehat{V}_R$, by setting the values at DoFs of $V_R$ that are shared between cells to be the average of the values from the neighbouring cells of the field in $\widehat{V}_R$.
The recovery operator is then expressed as
% \begin{linenomath}
\begin{equation}
\mathcal{R}=\mathcal{A}\widehat{\mathcal{P}}_R.
\end{equation}
% \end{linenomath}
As $V_R\nsubseteq V_H$, in place of the injection operator $\mathcal{I}_H$ we simply use a Galerkin projection $\mathcal{P}_H : V\to V_H$.
Then the whole reconstruction operator $\mathcal{J}$ can be expressed as
% \begin{linenomath}
\begin{equation}
\mathcal{J} := \mathcal{I}_H+\mathcal{P}_H\mathcal{R} - \mathcal{I}_H\mathcal{P}_L\mathcal{P}_H\mathcal{R},
\end{equation}
% \end{linenomath}
Again, the addition of $\mathcal{I}_H-\mathcal{I}_H\mathcal{P}_L\mathcal{P}_H\mathcal{R}$ ensures that the whole operation will be reversible in the absence of transport, as
% \begin{linenomath}
\begin{align*}
\mathcal{P}_L\mathcal{J} & =
\mathcal{P}_L\mathcal{I}_H + \mathcal{P}_L\mathcal{P}_H\mathcal{R}
- \mathcal{P}_L\mathcal{I}_H\mathcal{P}_L\mathcal{P}_H\mathcal{R} \\
& = \mathcal{P}_L\mathcal{I}_H + \mathcal{P}_L\mathcal{P}_H\mathcal{R}
- \mathcal{P}_L\mathcal{P}_H\mathcal{R} \\
& = \mathcal{P}_L\mathcal{I}_H ,
\end{align*}
% \end{linenomath}
which when acting upon a field in $V_L$ is the identity operator, since $V_L\subset V_H$.

\subsection{Choice of function spaces} \label{sec:space choices}
Armed with the extension of the recovery scheme presented in Section \ref{sec:recovery new scheme}, now consider the motivating case when $V_L$ is the lowest-order $H(\mathrm{div})$ Raviart-Thomas space for quadrilateral cells.
In general, there will be multiple possible choices for $V_R$ and $V_H$ that will satisfy the requirements presented in Section \ref{sec:recovery new scheme}, but here we only present the specific choices that are demonstrated in Section \ref{sec:results}.
These spaces are illustrated in Table \ref{tab:recovered spaces}.\\
\\
The general strategy is to choose $V_H$ to have the same continuity properties as $V_L$, but with increased polynomial order.
For the $H(\mathrm{div})$ Raviart-Thomas spaces, the components of the vector field that are normal to cell edges are continuous, and these already have a higher-order representation.
However the components that are tangential to cell edges are discontinuous with a lower-order representation.
Therefore we choose $V_R$ to be the higher-order $H(\mathrm{curl})$ space corresponding to $V_H$, whose tangential components are continuous between cells. \\
\\
On quadrilateral cells, $V_L$ is $\mathrm{RTc}^f_1$, which has a single DoF for each edge of the cell.
The higher-order space $V_H$ is $\mathrm{RTc}^f_2$, so from the same family as $V_L$ but with higher polynomial order.
The recovered space $V_R$ is the $H(\mathrm{curl})$ form, $\mathrm{RTc}^e_2$. \\
\\
For the transport operator $\mathcal{T}$, this work uses the benchmark upwind discretisation \eqref{eqn:dg upwind} in the higher-order space $V_H$, combined with a trapezoidal time discretisation that will be described in Section \ref{sec:results}.

\begin{table}[h!]
\centering
\begin{tabular}{m{0.2\textwidth}|m{0.2\textwidth}|m{0.2\textwidth}}
  \multicolumn{1}{c|}{$V_L: \mathrm{RTc}^f_1$}
& \multicolumn{1}{c|}{$V_R: \mathrm{RTc}^e_2$}
& \multicolumn{1}{c}{$V_H: \mathrm{RTc}^f_2$}  \\
\hline
  \includegraphics[width=0.2\textwidth]{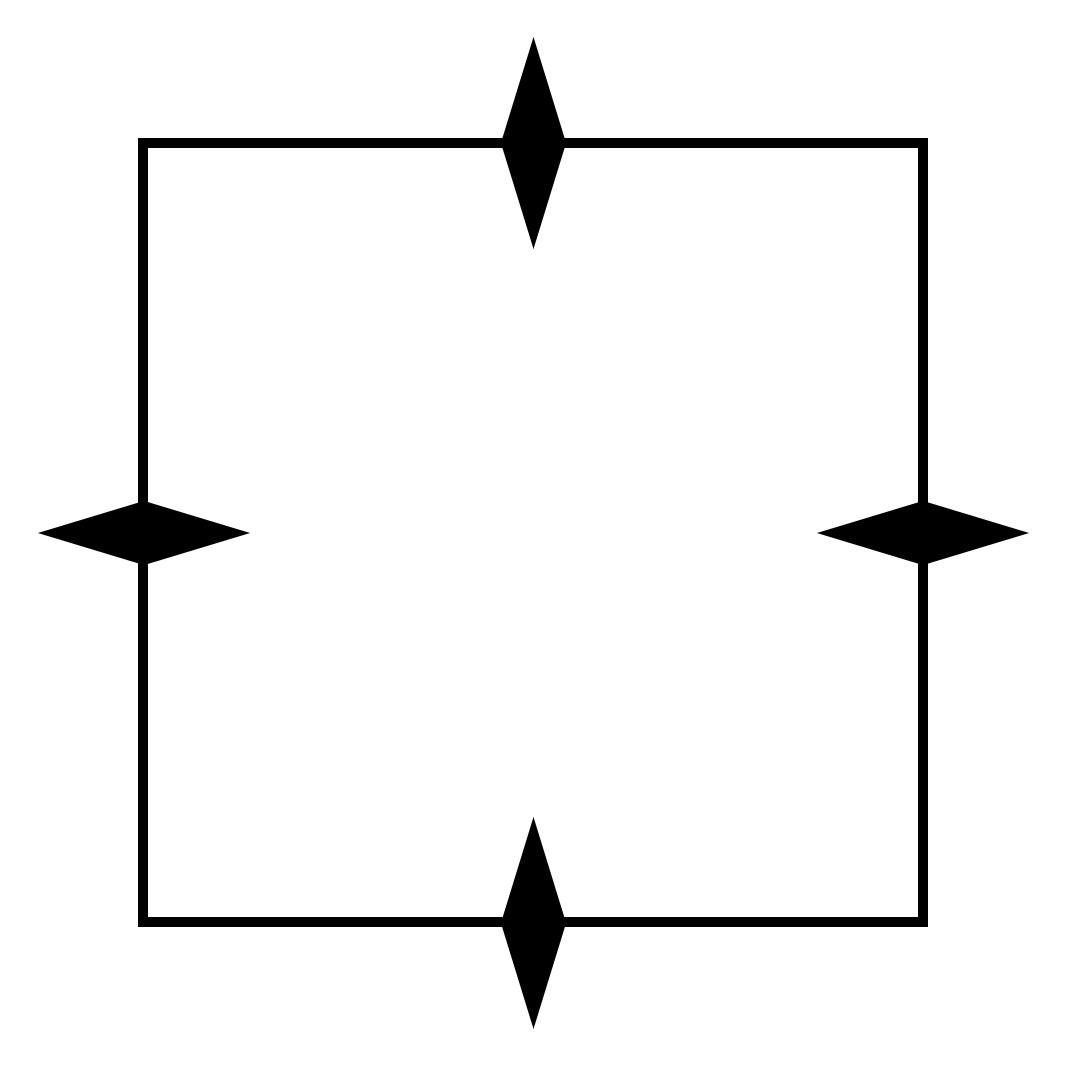}
& \includegraphics[width=0.2\textwidth]{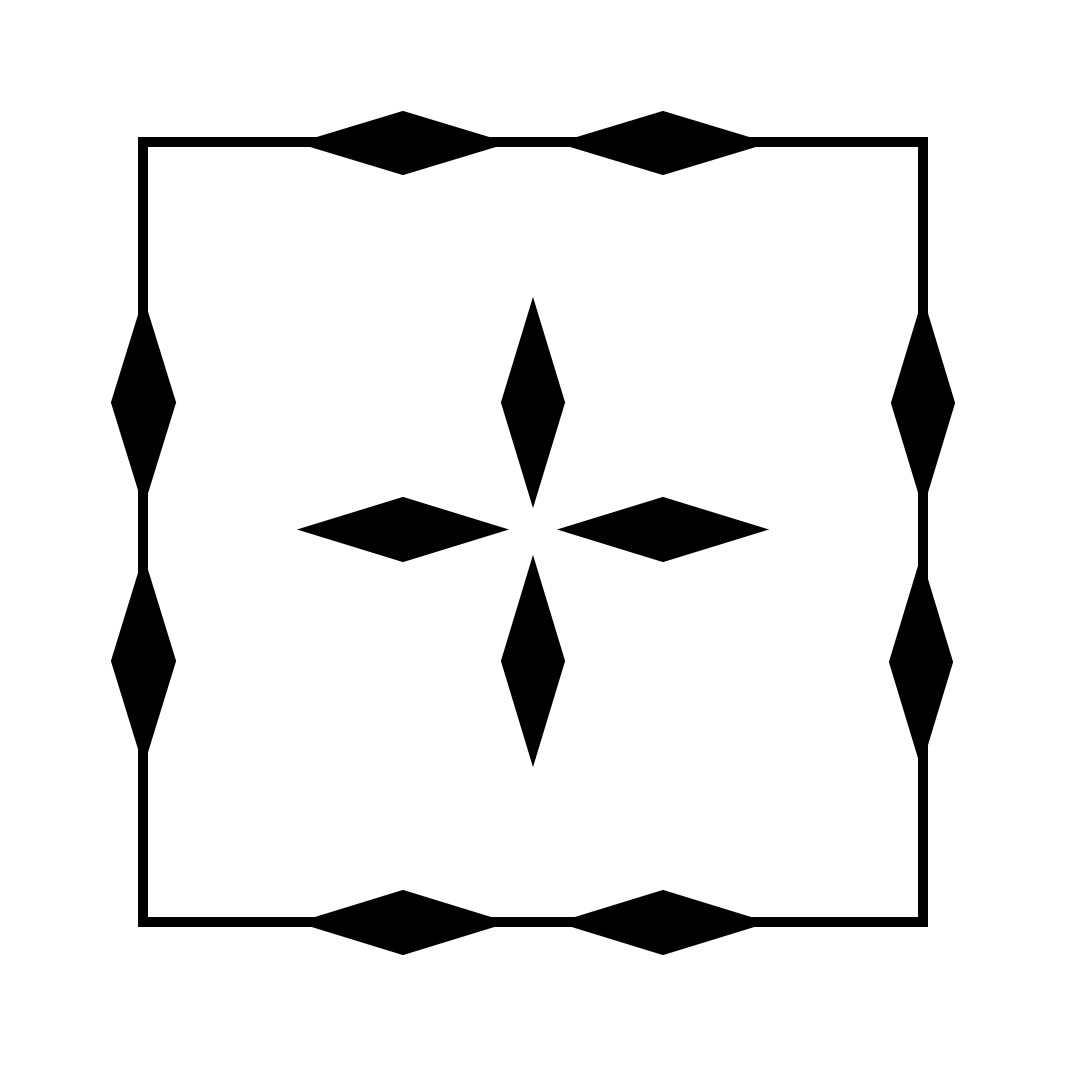}
& \includegraphics[width=0.2\textwidth]{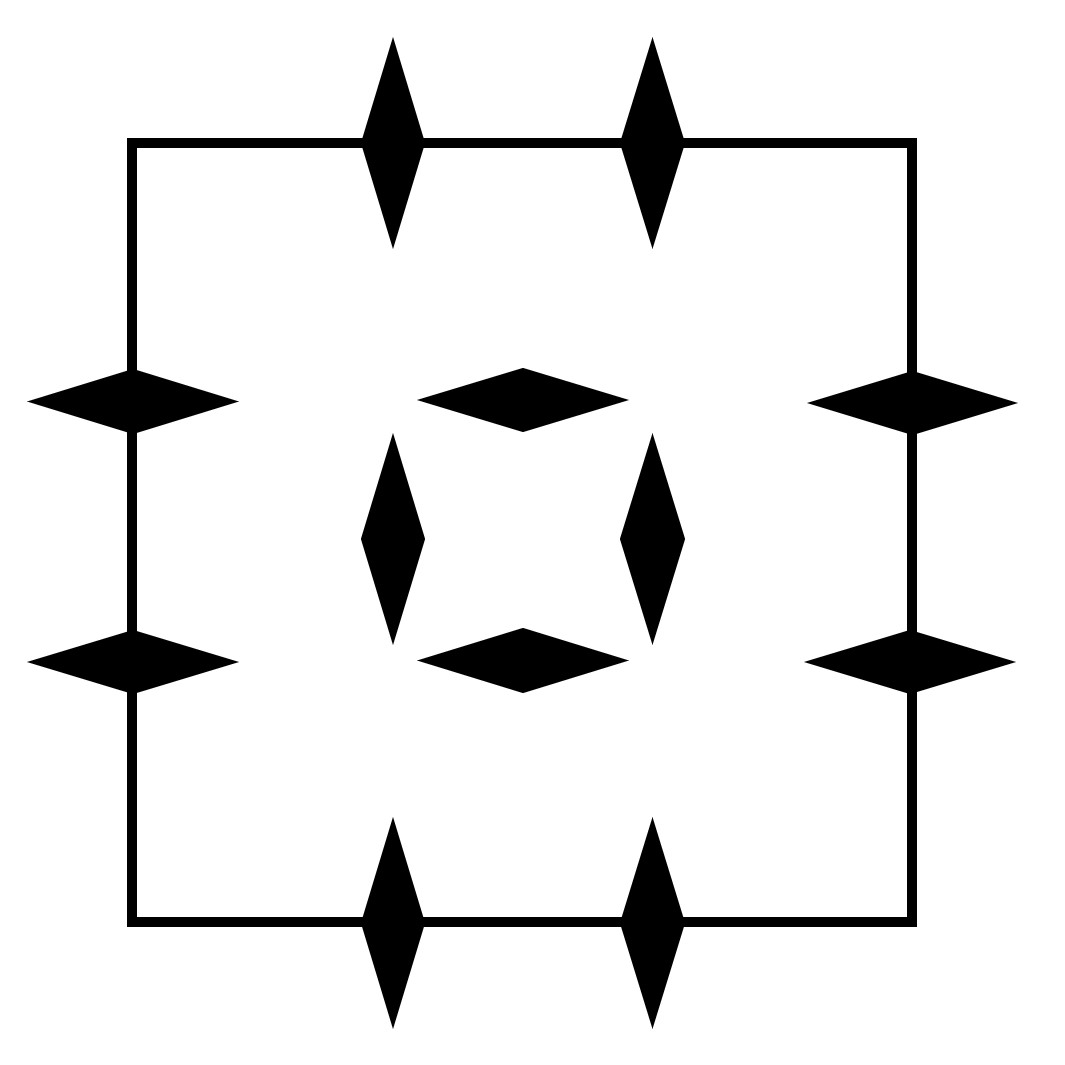} \\
    \end{tabular}
    \caption{Representations of the finite elements discussed in Section \ref{sec:space choices} for use in the recovered finite element method for transporting fields in the lowest-order Raviart-Thomas spaces.
    These are the specific choices of element that are used in the demonstrations of Section \ref{sec:results}.
    The diamonds represent the DoFs of the element, and whether the DoFs describe the components of the vector that are tangential or normal to the cell edges.
    }
    \label{tab:recovered spaces}
\end{table}

\section{Vorticity Form} \label{sec:vorticity}
As discussed at the end of Section \ref{sec:benchmark scheme}, the low order of accuracy of transport of fields in $\mathrm{RTc}^{f}_1$ can be attributed to the representation of the components of $\bm{F}$ that are perpendicular to the transporting velocity $\bm{v}$, as these components are only constant within a cell in the direction of $\bm{v}$.
When expressing the vector transport equation in the vorticity form \eqref{eqn:vec transport -- vorticity}, the transport of these components is captured in part by using the vorticity $\zeta$ through the $\zeta\bm{v}^\perp$ term.
This vorticity can be expressed weakly in a space $V_\zeta$ through
% \begin{linenomath}
\begin{equation} \label{eqn:zeta_def_discr}
\int_\varOmega \eta \; \zeta \dx{x} =
- \int_\varOmega \grad{^\perp\eta}\bm{\cdot F}\dx{x}, \qquad \forall\eta\in V_\zeta,
\end{equation}
% \end{linenomath}
where for simplicity terms associated with the boundary of the domain have been neglected\footnote{For a discussion including boundary terms, see \cite{bauer2018energy} and \cite{wimmer2020energy}.}.
In the compatible finite element framework, with $\bm{F}\in \mathrm{RTc}^f_k$, the space $V_\zeta$ is the continuous Galerkin space $\mathrm{CG}_k$, as for all $\eta\in\mathrm{CG}_k$ then $\bm{\nabla}^\perp\eta\in\mathrm{RTc}^f_k$.
Thus when using the lowest-order elements, $\zeta$ is piece-wise linear, which suggests that a formulation using the vorticity could improve the accuracy of the transport of $\bm{F}\in\mathrm{RTc}^f_1$. \\
\\
Due to its favourable properties with regards to the system's total energy and potential vorticity budgets \cite{ringler2010unified}, the vector-invariant form \eqref{eqn:vec transport -- vec inv} of the momentum equation is popular in NWP models (e.g.\ \cite{skamarock2012multiscale, zangl2015icon}).
It has also been used in the context of the compatible finite element 
discretisations of the shallow-water equations with vorticity as an auxiliary variable by
\cite{bauer2018energy}, \cite{wimmer2020energy}, \cite{mcrae2014energy}, and \cite{lee2021petrov}.
In \cite{mcrae2014energy}, a modified vorticity is diagnosed from the velocity field and used in its transport.
The modification to the vorticity improves the stability of the vector transport term by dissipating enstrophy without compromising on energy conservation.
This is known as the anticipated potential vorticity method (APVM).
However, this approach is not consistent with the continuous equations, in the sense that strong solutions do not necessarily satisfy the discrete equations. To overcome this problem, an extension to the APVM term was developed by \cite{bauer2018energy}, which used a mixed finite element problem to solve simultaneously for the velocity and the vorticity evolution equations, where the latter is stabilised by an SUPG method. For a recent comparison including the APVM and SUPG methods, see \cite{lee2022comparison}.
Finally, \cite{wimmer2020energy} introduced a methodology to apply an SUPG-based stabilisation method to more general vorticity evolution equations, which may contain additional terms such as ones arising from temperature gradients.
The methodology relies on the SUPG method's residual-based form.
In this section, we present a mixed velocity-vorticity approach for the vector transport equation similar to the setup of \cite{bauer2018energy}.
However, here we follow a ``black box'' approach in which the non-transport terms in the broader equation set are not included in the vorticity transport equation. To achieve this, we apply a residual-based setup akin to the one presented in \cite{wimmer2020energy}.
\\
\\
In order to derive an evolution equation for the vorticity, we consider the vorticity form of the vector transport equation
% \begin{linenomath}
\begin{equation}\label{eqn:vec transport -- vorticity sec_vort}
\pfrac{\bm{F}}{t} + \zeta \bm{v}^\perp
+ \frac{1}{2}\grad{(\bm{v\cdot F})}
+ \bm{G}(\bm{F}) = \bm{0},
\end{equation}
% \end{linenomath}
with
% \begin{linenomath}
\begin{equation}
    \bm{G}(\bm{F}) = \frac{1}{2}\left[\left(\bm{\nabla} \bm{F}\right)\bm{\cdot v}
-\left(\bm{\nabla} \bm{v}\right)\bm{\cdot F}  \right]. \label{eqn:G_cts}
\end{equation}
% \end{linenomath}
Applying the $\bm{\nabla}^\perp\bm{\cdot}$ operator and using $\bm{\nabla^\perp \cdot \nabla}f = 0$ for scalar $f$ yields
a vorticity equation of the form
% \begin{linenomath}
\begin{equation}\label{eqn:vort transport}
\pfrac{\zeta}{t} + \bm{\nabla} \bm{\cdot} (\zeta \bm{v})
+ \bm{\nabla^\perp \cdot} \bm{G}(\bm{F}) = \bm{0},
\end{equation}
% \end{linenomath}
noting that we applied the identity $\bm{a}^\perp \bm{\cdot b}^\perp = \bm{a} \bm{\cdot b}$, for any vectors $\bm{a}, \bm{b}$, to obtain
\begin{equation}
    \bm{\nabla}^\perp \bm{\cdot} (\zeta \bm{v}^\perp) = \bm{\nabla} \bm{\cdot} (\zeta \bm{v}). \label{cancel_perps}
\end{equation}
We arrive at our discretisation by multiplying \eqref{eqn:vec transport -- vorticity sec_vort} and \eqref{eqn:vort transport} by test functions $\bm{\gamma}\in V_F$ and $\eta\in V_\zeta$, which yields a mixed finite element problem with two equations to be solved simultaneously:
\begingroup
\addtolength{\jot}{2mm}
% \begin{linenomath}
\begin{subequations} \label{eqn:F_zeta_discr_nostab}
\begin{align}
& \int_\varOmega \bm{\gamma} \bm{\cdot} \pfrac{\bm{F}}{t} \dx{x} +
\int_\varOmega \bm{\gamma} \bm{\cdot} (\zeta \bm{v}^\perp)\dx{x} - 
\frac{1}{2} \int_\varOmega (\bm{v\cdot F}) \; (\bm{\nabla \cdot} \bm{\gamma})\dx{x}
+ \bm{G}'(\bm{F}; \bm{\gamma}) = 0, &\forall \bm{\gamma} \in V_F, \label{eqn:F_discr_nostab} \\
& \int_\varOmega \eta \pfrac{\zeta}{t} \dx{x} - \int_\varOmega \bm{\nabla}\eta \bm{\cdot} (\zeta \bm{v})\dx{x} - \bm{G}'(\bm{F}; \bm{\nabla^\perp} \eta) = 0, & \forall \eta \in V_\zeta, \label{eqn:zeta_discr_nostab}
\end{align}
\end{subequations}
% \end{linenomath}
\endgroup
where the initial discrete vorticity $\zeta$ is defined by \eqref{eqn:zeta_def_discr}, and $\bm{G}'$ is a weak discretisation of $\bm{G}$, whose specific value is postponed to later in this section. Note that to arrive at the above weak vorticity equation, we applied integration by parts according to
% \begin{linenomath}
\begin{align}
    \int_\Omega \eta \bm{\nabla \cdot} (\zeta \bm{v}) \dx{x} = - \int_\Omega \bm{\nabla} \eta \bm{\cdot} (\zeta \bm{v}) \dx{x} && \forall \eta \in V_\zeta,
\end{align}
% \end{linenomath}
which does not include any additional facet integral terms since the choice of finite element spaces ensures that the normal component of $\zeta \bm{v}$ is continuous.
Since \eqref{eqn:F_zeta_discr_nostab} is true for all $\bm{\gamma}$, it is also true for $\bm{\gamma}=-\grad{^\perp}\eta$ which recovers \eqref{eqn:zeta_discr_nostab} (by cancelling perpendicular operations similar to \eqref{cancel_perps}). This means that by solving these equations simultaneously, the evolution of the discrete $\bm{F}$ and its vorticity are kept consistent.
\\
\\
In the form of \eqref{eqn:zeta_discr_nostab}, the discrete vorticity evolution equation does not contain any transport stabilisation measures, and we may therefore expect it to be vulnerable to grid-scale oscillations.
In the context of a shallow-water model this can correspond to a lack of dissipation of enstrophy, which naturally cascades to fine scales but gets trapped at the grid scale without a mechanism to dissipate it \cite{mcrae2014energy}.
Since the vorticity is discretised as a $\mathrm{CG}_k$ field, this can be remedied by using a stabilisation based on the SUPG method.
\\
\\
The usual Petrov-Galerkin approach to applying an SUPG stabilisation is to adjust the test function to include a transport contribution via
% \begin{linenomath}
\begin{equation}
\eta \; \to \; \eta + \tau \bm{v} \bm{\cdot \nabla} \eta, \label{standard_supg}
\end{equation}
% \end{linenomath}
where $\tau$ denotes a suitable stabilisation parameter with dimensions of time. However, modifying only the test function for \eqref{eqn:zeta_discr_nostab} breaks the consistency between the evolution equations of $\bm{F}$ and $\zeta$.
Instead, we
consider a residual-based approach like those used by \cite{wimmer2020energy}.
This uses the residual of the strong form of the vorticity equation:
% \begin{linenomath}
\begin{equation}
    \zeta_{res} = \pfrac{\zeta}{t} + \bm{\nabla} \bm{\cdot} (\zeta \bm{v})
+ \bm{\nabla^\perp \cdot} \bm{G}(\bm{F}).
\end{equation}
% \end{linenomath}
Then, the vorticity appearing in the discretisation \eqref{eqn:F_zeta_discr_nostab} is modified to give
\begingroup
\addtolength{\jot}{2mm}
% \begin{linenomath}
\begin{subequations} \label{eqn:F_zeta_discr}
\begin{align}
& \int_\varOmega \bm{\gamma} \bm{\cdot} \pfrac{\bm{F}}{t} \dx{x} +
\int_\varOmega \bm{\gamma} \bm{\cdot} (\zeta^\ast \bm{v}^\perp)\dx{x} - 
\frac{1}{2} \int_\varOmega (\bm{v\cdot F}) \; (\bm{\nabla \cdot} \bm{\gamma})\dx{x}
+ \bm{G}'(\bm{F}; \bm{\gamma}) = 0, &\forall \bm{\gamma} \in V_F, \label{eqn:F_discr} \\
& \int_\varOmega \eta \pfrac{\zeta}{t} \dx{x} - \int_\varOmega \bm{\nabla} \eta \bm{\cdot} (\zeta^\ast \bm{v})\dx{x} - \bm{G}'(\bm{F}; \bm{\nabla^\perp} \eta) = 0, & \forall \eta \in V_\zeta, \label{eqn:zeta_discr}
\end{align}
\end{subequations}
% \end{linenomath}
\endgroup
with $\zeta^\ast = \zeta - \tau \zeta_{res}$.
Note that after discretisation, the differential operations occurring in $\zeta_{res}$ are applied cell-wise, including those of $\bm{G}(\bm{F})$ as defined by \eqref{eqn:G_cts}. There is then choice in the time discretisation; this is discussed briefly in Section \ref{sec:results}.
\\
\\ 
We conclude the description of the SUPG stabilisation with the following four observations.
First, the choice of residual $\zeta_{res}$ ensures that the discretisation \eqref{eqn:F_zeta_discr} is consistent with the strong equation of the vorticity evolution, as then $\zeta_{res} = 0$ and $\zeta^\ast$ reduces to $\zeta$.
At the same time the evolution equations for $\bm{F}$ and $\zeta$ are still consistent with one another.
Secondly, the modification to the vorticity used in \eqref{eqn:F_zeta_discr} has a stabilising effect akin to the more standard SUPG modification \eqref{standard_supg}.
This can be seen by setting $\eta = \zeta$ in \eqref{eqn:zeta_discr}, which leads to a non-positive definite term of the form
% \begin{linenomath}
\begin{equation}
    \frac{1}{2}\frac{d}{dt}\|\zeta\|_2^2 = \int_\varOmega \zeta \pfrac{\zeta}{t}\dx{x} = \cdots - \|\sqrt{\tau} \bm{v} \cdot \bm{\nabla} \zeta \|_2^2, \label{vort_dissipation}
\end{equation}
% \end{linenomath}
on the equation's right-hand side, showing that as expected for the SUPG method, there is potential for vorticity dissipation along the direction of the flow.
Thirdly, whilst the above formulation allows for the dissipation of vorticity, it does not necessarily dissipate the divergence field. If the latter field is large, additional stabilisation mechanisms may be required for the transport of $\bm{F}$.
An example for this would be an interior penalty term \cite{burman2004edge}, based on the divergence field $\bm{\nabla \cdot \bm{F}} \in \text{DG}_{k-1}$.
For the type of shallow-water scenarios typically considered in numerical weather prediction, the divergence field is small, and no such additional mechanism is required.
Lastly, if the term $\bm{G}$ and its weak discrete version $\bm{G}'$ are equal to zero -- as will be the case if the advecting velocity $\bm{v}$ is set equal to $\bm{F}$ -- then the vorticity evolution equation can be rewritten in standard SUPG form
% \begin{linenomath}
\begin{align}
    \int_\varOmega \left( \eta + \tau \bm{v} \cdot \bm{\nabla} \eta \right) \left( \pfrac{\zeta}{t} + \bm{\nabla} \bm{\cdot} (\zeta \bm{v}) \right)\dx{x}  = 0, && \forall \eta \in V_\zeta.
\end{align}
% \end{linenomath}
Note that to arrive at the above equation, we applied integration by parts, which does not lead to any additional facet integrals as mentioned above when deriving \eqref{eqn:zeta_discr_nostab}. When $\bm{G}'$ is non-zero, the non-equivalence of the Petrov-Galerkin and residual-based approaches is a necessity arising from formulating \eqref{eqn:zeta_discr} in a manner consistent with \eqref{eqn:F_discr}.
This non-equivalence can also be found in other applications in the literature, such as SUPG discretisations of the Navier-Stokes equations.
In the latter case, a residual-based formulation may be preferred in order to avoid a double-derivative applied to the SUPG-modified test function in the weak diffusion term \cite{elman1996iterative}.\\
\\
It remains to describe the weak, discrete operator $\bm{G}'$.
In order to stabilise the gradient terms occurring in $\bm{G}$, an upwind formulation is used, so that for test functions $\bm{w}$
% \begin{linenomath}
\begin{equation}
\bm{G}'(\bm{F}; \bm{w}) = \frac{1}{2} \int_\varOmega \left(
\bm{F\cdot}\left[\grad{}\bm{\cdot}(\bm{v}\otimes \bm{w})\right] -
\bm{v \cdot}\left[\grad{}\bm{\cdot}(\bm{F} \otimes \bm{w})\right]
\right)\dx{x}
+\frac{1}{2}\int_\Gamma \left(\bm{w}^+\bm{\cdot}\widehat{\bm{n}}^+\right)
\left(\llbracket \bm{v}\rrbracket_+ \bm{\cdot}\bm{F}^\dagger - 
\llbracket \bm{F}\rrbracket_+ \bm{\cdot}\bm{v}^\dagger\right) \dx{S}.
\end{equation}
% \end{linenomath}
Again, terms associated with the boundaries of the domain have been neglected.
As with the benchmark \eqref{eqn:dg upwind}, a correction could be added to project the upwind term into the tangent bundle. \\
\\
Finally, it should be stressed that in the context of a shallow-water model, alternative variables to the (relative) vorticity $\zeta$ are the \textit{absolute vorticity} $\omega=\bm{\nabla}^\perp\bm{\cdot u} +f$ or \textit{potential vorticity} $q = (\bm{\nabla^\perp \cdot u} + f)/h$.
As the potential vorticity is conserved along the flow, it is often preferred to $\zeta$, as in the case of the compatible finite element discretisations in \cite{mcrae2014energy} and \cite{bauer2018energy}.
In particular, the APVM and SUPG stabilisations derived in the aforementioned papers dissipate the enstrophy $hq^2$, while conserving the system's total energy. These works also solved the vorticity evolution equation corresponding to the whole shallow-water equation for the velocity \eqref{eqn:shallow-water-momentum}, whereas in this section we consider only the transport part; as mentioned in Section \ref{sec:shallow-water background} the motivation here is to find a ``black box'' to solve the vector transport equation.
The addition of the SUPG stabilisation is still consistent \textit{within the transport step}, and any errors arising from not applying SUPG to the whole equation will do so in the form of a splitting error in time.
Note that this SUPG setup is different to the ones used in \cite{bauer2018energy} and \cite{wimmer2020energy}.
In the former, a different vorticity variable is used, leading to a forcing contribution of the form $g\grad{(h+h_b)}$, which vanishes in the vorticity evolution equation (since $\bm{\nabla^\perp \cdot} \grad{(h+h_b)} = 0$).
In the latter, there is no time splitting, and a forcing contribution $\bm{\nabla^\perp \cdot} \bm{J}$ for some baroclinic forcing terms $\bm{J}$ is included in the vorticity equation's residual. While these approaches avoid errors due to time splitting and lead to additional conservation properties such as energy conservation, they require additional information from the equations and cannot be used as ``black box'' vector transport methods.
In particular, this is a drawback for code implementations: while a ``black box'' setup can be used for a variety of different equation sets, in the specific setups of \cite{bauer2018energy} and \cite{wimmer2020energy}, the transport implementation has to be adjusted each time the overall equation sets are changed.

\section{Numerical Results} \label{sec:results}
This section demonstrates the schemes presented in Sections \ref{sec:recovery} and \ref{sec:vorticity}, through some transport-only tests in Section \ref{sec:transport only} and in the context of a shallow-water model in Section \ref{sec:shallow-water}.
The new schemes are compared with the benchmark scheme of Section \ref{sec:benchmark scheme}, all applied to the lowest-order quadrilateral Raviart-Thomas elements. \\
\\
Throughout this section, the equations are discretised in time using the trapezoidal rule.
If the discretisation of the integrated transport term is given by $\mathcal{G}\left[\bm{\gamma},\bm{v},\bm{F}\right]$, then the value of $\bm{F}$ at the $(n+1)$-th time step is found from
% \begin{linenomath}
\begin{equation} \label{eqn:trapezoidal}
\int\bm{\gamma\cdot}\left(\bm{F}^{n+1} -\bm{F}^{n}\right)\dx{x} = \frac{\Delta t}{2}\left(\mathcal{G}\left[\bm{\gamma},\bm{v},\bm{F}^n\right] + \mathcal{G}\left[\bm{\gamma},\bm{v},\bm{F}^{n+1}\right] \right), \qquad \forall \bm{\gamma}\in V_F.
\end{equation}
% \end{linenomath}
This yields a matrix-vector problem for $\bm{F}^{n+1}$ which is then solved to obtain the transported solution.
For the mixed vorticity scheme of Section \ref{sec:vorticity}, the integrated transport terms for $\bm{F}$ and $\zeta$ are given respectively by $\mathcal{G}\left[\bm{\gamma},\bm{v},\bm{F},\zeta\right]$ and $\mathcal{H}\left[\eta,\bm{v},\bm{F},\zeta\right]$, so that the trapezoidal rule is given by
% \begin{linenomath}
\begin{subequations}
\begin{align}
\int\bm{\gamma\cdot}\left(\bm{F}^{n+1} -\bm{F}^{n}\right)\dx{x} = \frac{\Delta t}{2}\left(\mathcal{G}\left[\bm{\gamma},\bm{v},\bm{F}^n,\zeta^{n}\right] + \mathcal{G}\left[\bm{\gamma},\bm{v},\bm{F}^{n+1},\zeta^{n+1}\right] \right), \qquad &\forall \bm{\gamma}\in V_F, \\
\int\eta\left(\zeta^{n+1} -\zeta^{n}\right)\dx{x} = \frac{\Delta t}{2}\left(\mathcal{H}\left[\eta,\bm{v},\bm{F}^n,\zeta^n\right] + \mathcal{H}\left[\eta,\bm{v},\bm{F}^{n+1},\zeta^{n+1}\right] \right), \qquad &\forall \eta\in V_\zeta.
\end{align}
\end{subequations}
% \end{linenomath}
To implement these schemes, we used the Firedrake software, \cite{Rathgeber2016}, which is a library for solving PDEs using finite element methods and is built on the PETSc solver library \cite{petsc-user-ref}.
Firedrake constructs the quadrilateral Raviart-Thomas elements as tensor-product elements \cite{McRae2016} and provides support for the hybridised solver \cite{gibson2020slate} used in the shallow-water model of Section \ref{sec:shallow-water}. The orthographic projections were plotted using the Cartopy python package \cite{Cartopy}. Finally, for the SUPG method used in the stabilised voriticity discretisation, we consider a stabilisation parameter of the form
% \begin{linenomath}
\begin{equation}
    \tau = \left(\lambda \frac{2}{\Delta t} + \frac{2|\bm{u}|}{\Delta x} \right)^{-1},
\end{equation}
% \end{linenomath}
for local mesh size $\Delta x$, and a tuning parameter $\lambda \ge 0$. The latter parameter can be seen to adjust the stabilisation's ``aggressiveness'' and in this section, we took $\lambda=0.5$; for details, see \cite{wimmer2020energy}.

\subsection{Transport-only tests} \label{sec:transport only}
Although the literature on numerical weather prediction contains many test cases for the transport of scalar fields, there are few for the transport of vector fields. This section describes two test cases on curved manifolds that may be used for assessing the convergence properties of transport schemes for vector fields.
\subsubsection{Deformation on the cylinder} \label{sec:cylinder test}
\begin{figure}[h!]
    \centering
    \includegraphics[width=0.95\textwidth]{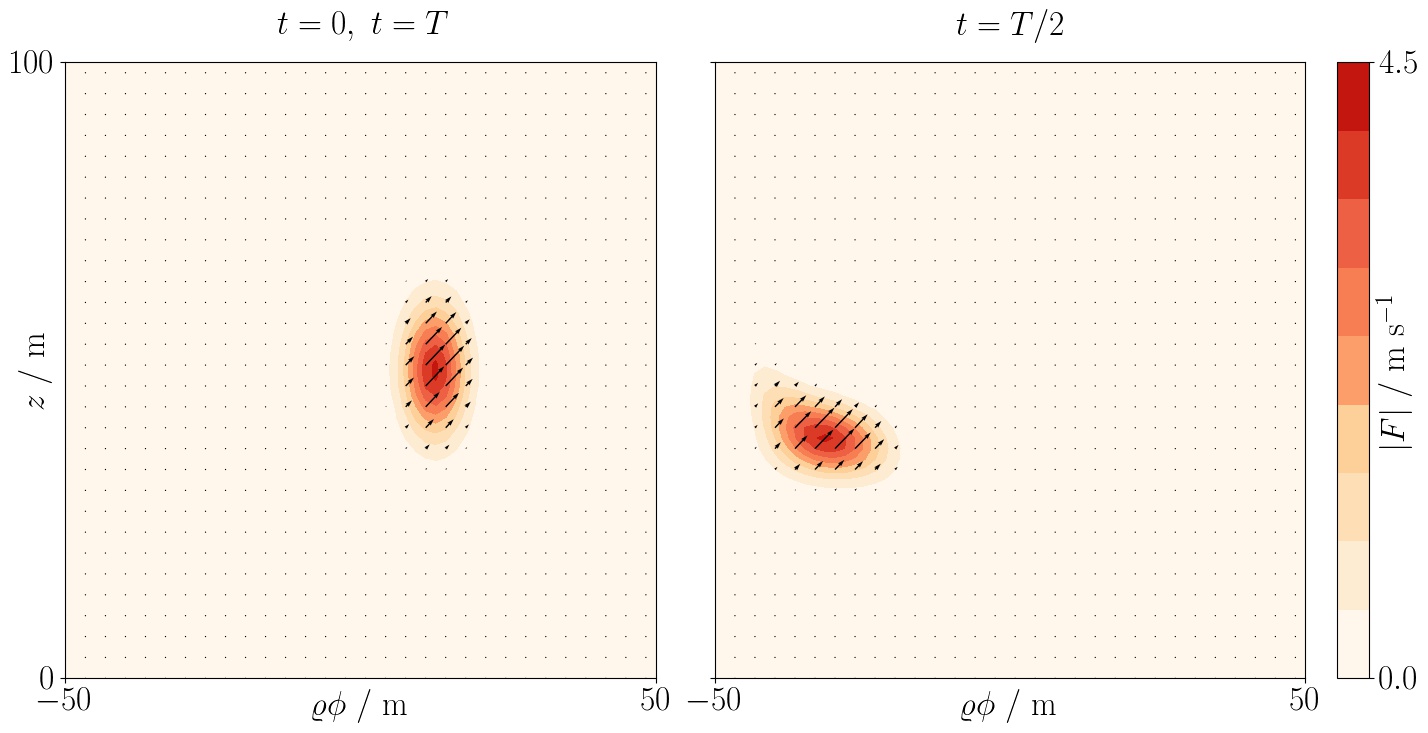}
    \caption{The transported field $\bm{F}$ in the deformational cylindrical transport test of Section \ref{sec:cylinder test}.
    The contours show the magnitude of $\bm{F}$, with the arrows indicating its direction.
    (Left) the initial condition, and true solution at $t=T$.
    (Right) a numerical computation of the deformed field at $t=T/2$.
    The contours are spaced at 0.5 m s$^{-1}$.
    }
    \label{fig:cylinder ics}
\end{figure}

\begin{figure}[h!]
    \centering
    \includegraphics[width=0.95\textwidth]{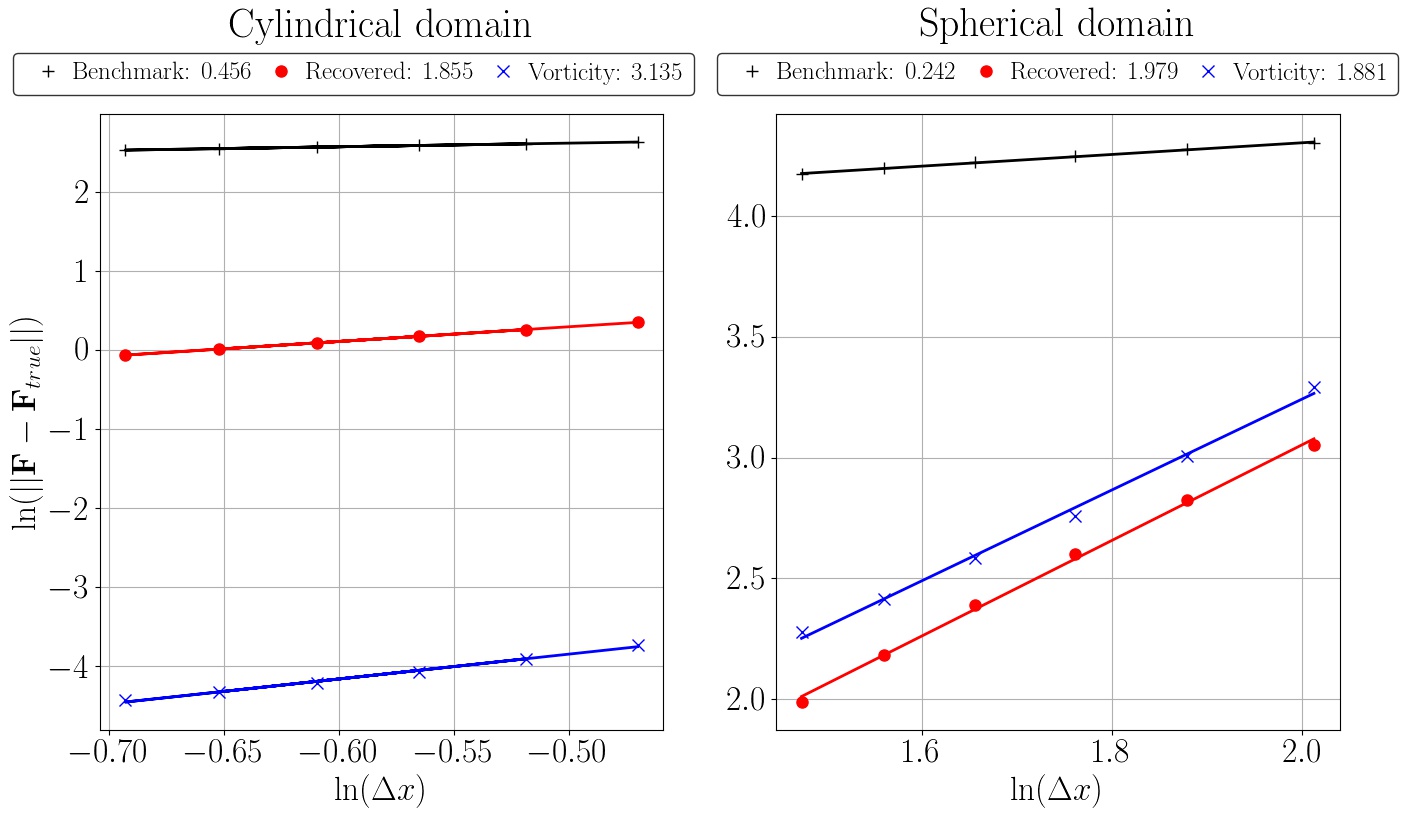}
    \caption{Convergence results for the transport test of Sections \ref{sec:cylinder test} and \ref{sec:sphere test}.
    The $L^2$ errors in the transported $\bm{F}$ field are plotted for a range of spatial resolutions.
    The benchmark case of Section \ref{sec:benchmark scheme} is compared with the recovered scheme of Section \ref{sec:recovery} and the vorticity scheme of \ref{sec:vorticity} with the SUPG stabilisation.
    The legends indicate the gradients of lines of best fit through the error measurements, which approximate the rate of convergence of the scheme.
    (Left) results for the cylindrical test, and (right) results for the spherical test.
    For both tests, the two new schemes demonstrate much better convergence than the benchmark case.
    }
    \label{fig:convergence transport}
\end{figure}

The surface of a cylinder is a curved manifold on which the vector transport equation does not have metric terms.
It is therefore straightforward to adapt existing transport tests to the cylinder using the standard format for convergence tests, in which the true final solution of a transported vector is equal to its initial state.
If the azimuthal and height coordinates are $\bm{x}=(\phi,z)$ and the radius of the cylinder is $\varrho$, the vector transport equation \eqref{eqn:vec transport -- advective} can be expressed in components as
% \begin{linenomath}
\begin{subequations}
\begin{align}
& \pfrac{F_\phi}{t} + \frac{v_\phi}{\varrho}\pfrac{F_\phi}{\phi} +
v_z\pfrac{F_z}{z} = 0, \\
& \pfrac{F_z}{t} + \frac{v_\phi}{\varrho}\pfrac{F_\phi}{\phi} +
v_z\pfrac{F_z}{z}  = 0.
\end{align}
\end{subequations}
% \end{linenomath}
Here the transporting velocity is inspired by the time-varying and deformational divergence-free flows from \cite{nair2010class} and \cite{lauritzen2012standard}, but adapted to the cylinder.
The cylinder has radius $\varrho$ and length $L$, which is periodic in the $z$ direction.
The time $t$ runs from 0 to $T$.
With speeds $U=2\pi\varrho/T$ and $W$, and a modified coordinate $\phi'=\phi-Ut/\varrho$, the transporting velocity is given by
% \begin{linenomath}
\begin{subequations}
\begin{align}
v_\phi & = U + 2\pi W\sin\left(\phi'\right)\sin\left(\frac{2\pi z}{L}\right)\cos\left(\frac{\pi t}{T}\right), \\
v_z & = \frac{W L}{\varrho}\cos\left(\phi'\right)\cos\left(\frac{2\pi z}{L}\right)\cos\left(\frac{\pi t}{T}\right).
\end{align}
\end{subequations}
% \end{linenomath}
With this flow the true solution at $t=T$ is equal to the initial condition.
As in \cite{nair2010class} and \cite{lauritzen2012standard}, the flow has a translational component to avoid fortuitous cancellation of errors.
The amount of deformation can be controlled by changing $W$ relative to $U$.
This flow can also be expressed using a stream function, but this will contain a jump on the periodic cylinder due to the translational component of the flow.
For our test we took $L=100$ m, $\varrho = L/(2\pi)$, $T = 100$ s and $W=U/10$.\\
\\
To describe the initial conditions, let the distance from a specific point $(\phi_c, z_c)$ on the cylindrical surface be defined via
% \begin{linenomath}
\begin{equation}
\ell^2(\phi,z) = \left(\cos^{-1}\left[\cos\left(\phi - \phi_c\right)\right] \right)^2 + \left(\cos^{-1}\left[\cos\left(\frac{2\pi(z-z_c)}{L}\right)\right] \right)^2.
\end{equation}
% \end{linenomath}
The initial condition uses a vector whose cylindrical components are both a Gaussian hill of size $F_0$, width $\ell_0$ and centred on $(\phi_c,z_c)$, taking:
% \begin{linenomath}
\begin{equation}
\bm{F} = 
\left(\widehat{\bm{e}}_\phi + \widehat{\bm{e}}_z\right)F_0\exp\left(-\ell^2(\phi,z)/\ell^2_0\right),
\end{equation}
% \end{linenomath}
where $\phi_c=\pi/4$, $z_c=L/2$, $\ell_0=1/10$, $F_0=3$ m s$^{-1}$.
This initial condition and a numerical solution at $t=T/2$ are displayed in Figure \ref{fig:cylinder ics}.\\
 \\
To perform a convergence test, the $L^2$ error was computed for the numerical solution against the true solution at $t=T$, for a range of spatial resolutions.
The same time step $\Delta t = 0.002$ s was used for all simulations, and the meshes were constructed of uniform quadrilateral cells.
Results of the convergence test comparing the benchmark scheme of Section \ref{sec:benchmark scheme} to the new schemes are shown in the left of Figure \ref{fig:convergence transport}.
Both schemes show a very clear improvement from the benchmark scheme, with the recovered scheme approaching second-order accuracy and the vorticity scheme (which used the SUPG stabilisation) even achieving some super-convergence.

\subsubsection{Solid body rotations on sphere} \label{sec:sphere test}
\begin{figure}[h!]
    \centering
    \includegraphics[width=0.95\textwidth]{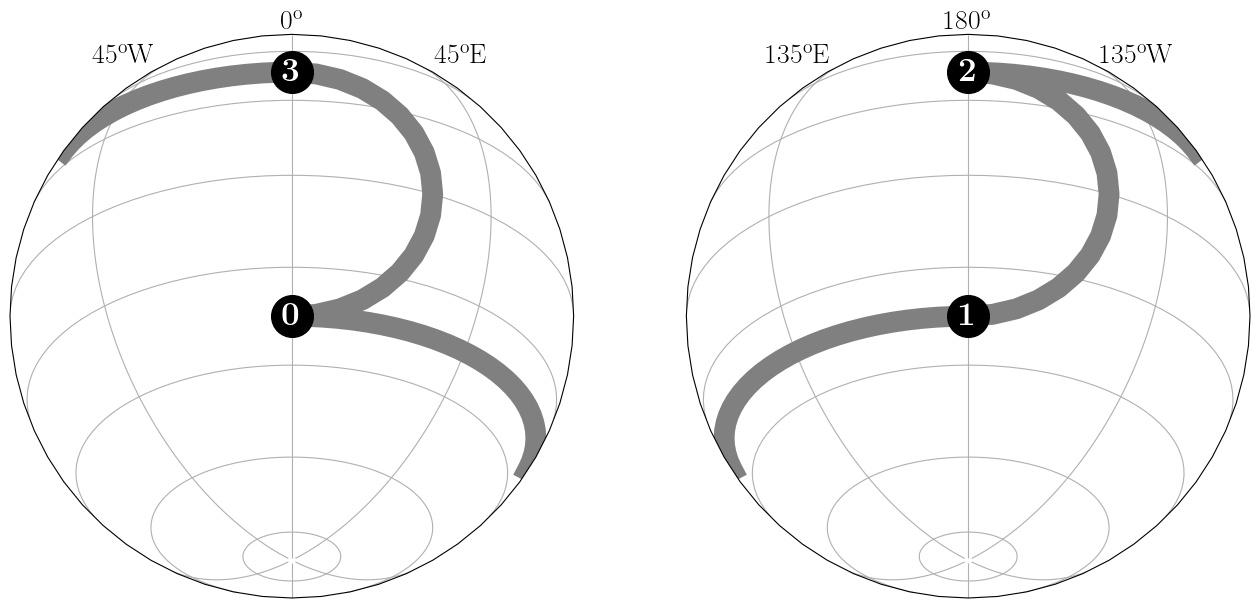}
    \caption{An illustration of the path taken by the transported field in the solid body rotation test presented in Section \ref{sec:sphere test}.
    The `front' of the sphere is shown on the left and the `back' on the right.
    The path, shown in grey, is broken into four stages, marked by the black circles.
    Each stage involves a solid body rotation: from points 0 to 1 and 2 to 3 this is a solid body rotation around the $z$-axis, while it is a solid body rotation around the $x$-axis from points 1 to 2 and 3 to 0.
    Taking this path, the effects induced by metric terms upon a transported vector cancel out, as any transport at a latitude $\vartheta$ is matched by equal transport at $-\vartheta$.
    Thus the true final solution is equal to the initial solution.}
    \label{fig:spherical path}
\end{figure}

\begin{figure}[h!]
\centering
\includegraphics[width=0.95\textwidth]{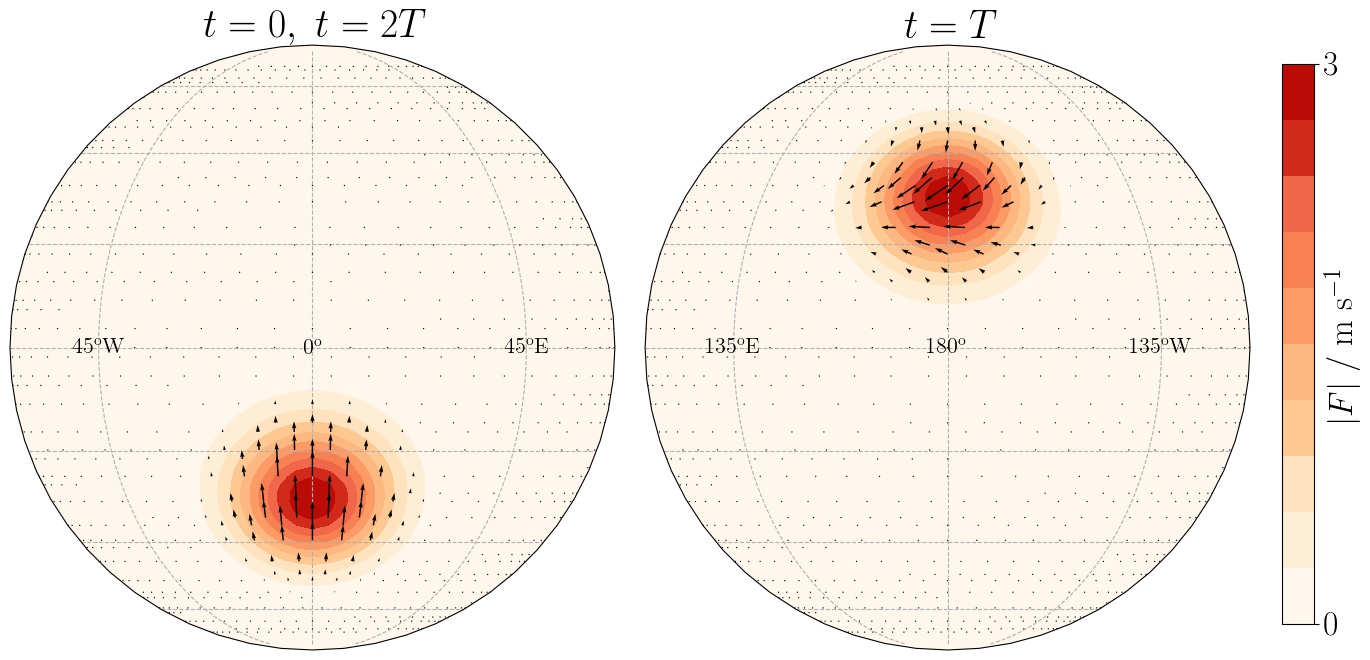}
\caption{The transported field $\bm{F}$ for the spherical transport test of Section \ref{sec:sphere test}.
The contours show the magnitude of $\bm{F}$, with the arrows indicating its direction.
(Left) the initial condition and true solution at $t=2T$, shown on the `front' of the sphere.
The meridional velocity is a Gaussian profile centred on $\lambda_c=0$ and $\vartheta_c=-\pi/6$.
(Right) a numerical computation of the transported field at $t=T$, shown on the `back' of the sphere.
We see the effect of the metric terms here on the direction of $\bm{F}$.
The contours are spaced at $0.3$ m s$^{-1}$.
}
\label{fig:spherical ics}
\end{figure}

Unlike on a cylindrical manifold, the vector transport equation does have metric terms on a spherical manifold.
If $(\lambda,\vartheta)$ are the longitude and latitude, and $r$ is the radius of the sphere, the advective form of the transport equation \eqref{eqn:vec transport -- advective} can be written as
% \begin{linenomath}
\begin{subequations} \label{eqn:spherical vec transport}
\begin{align}
& \pfrac{F_\lambda}{t} + \frac{v_\lambda}{r\cos\vartheta}\pfrac{F_\lambda}{\lambda} +
\frac{v_\vartheta}{r}\pfrac{F_\lambda}{\vartheta} - \frac{v_\lambda F_\vartheta \tan\vartheta}{r} = 0, \\
& \pfrac{F_\vartheta}{t} + \frac{v_\lambda}{r\cos\vartheta}\pfrac{F_\vartheta}{\lambda} +
\frac{v_\vartheta}{r}\pfrac{F_\vartheta}{\vartheta} + \frac{v_\lambda F_\lambda \tan\vartheta}{r} = 0.
\end{align}
\end{subequations}
% \end{linenomath}
The presence of the metric terms in these equations makes the design of a convergence test difficult, as any zonal component of $\bm{v}$ will cause the rotation of $\bm{F}$ at a rate depending on the latitude.
One strategy to avoid this is to explicitly add the metric terms as a forcing to the equation, but this requires a discretisation of the metric terms themselves which can confuse the interpretation of any results.
Another strategy is to use an exactly reversing flow to cancel out the effects of the metric terms, but this could also result in the fortuitous cancellation of dispersion errors.\\
\\
Here we present a spherical convergence test that avoids these issues by composing four solid body rotations to reverse the effects of the metric terms.
First, $\bm{F}$ is initialised with a smooth profile centred at $(\lambda_c, \vartheta_c)$, taking $\lambda_c=0$ and $\vartheta_c=-\pi/6$. 
Using the usual definition of distance on a spherical surface,
% \begin{linenomath}
\begin{equation}
\ell(\lambda,\vartheta) = \cos^{-1}\left[\sin\vartheta_c\sin\vartheta+\cos\vartheta_c\cos\lambda_c\cos\lambda+\cos\vartheta_c\sin\lambda_c\sin\lambda \right],
\end{equation}
% \end{linenomath}
the initial condition is
% \begin{linenomath}
\begin{equation}
F_\lambda = 0, \quad F_\vartheta = F_0 \exp\left(-\ell^2(\lambda,\vartheta)/\ell_0^2\right),
\end{equation}
% \end{linenomath}
with $F_0=3$ and $\ell_0=1/4$.
This is displayed in Figure \ref{fig:spherical ics}. \\
\\
The transporting velocity is made by composing four solid body rotations, each performing half of a rotation of the profile around an axis.
The first half-rotation is around the $z$-axis, leaving a profile that should be centred on $\lambda=\pi$.
Then the velocity is changed to perform a half-rotation around the $x$-axis, rotating the profile from the southern hemisphere to the northern hemisphere.
The third rotation uses the same winds as the first, rotating the profile again around the $z$-axis. 
By performing the same solid body rotation again, but this time with the profile in the northern hemisphere instead of the southern, the metric effects induced by the first half-rotation will be cancelled out.
Finally, another half-rotation is completed around the $x$-axis, which reverses the effects of the metric terms from the previous rotation around the $x$-axis.
The resulting path around the sphere is illustrated in Figure \ref{fig:spherical path}.\\
\\
The transporting velocity can be summarised in $(\lambda,\vartheta)$ components as
% \begin{linenomath}
\begin{equation} \label{eqn:rotations}
\begin{array}{lll}
v_\lambda = U\cos\vartheta, & v_\vartheta = 0,
& \mathrm{for} \ 0\leq t\leq T/2 \ \mathrm{and} \ T < t \leq 3T/2, \\
& & \\
v_\lambda = -U\cos\lambda\sin\vartheta,
& v_\vartheta = U\sin\lambda(\cos^2\vartheta-\sin^2\vartheta)
& \mathrm{for} \ T/2 < t \leq T \ \mathrm{and} \ 3T/2 < t \leq 2T.
\end{array}
\end{equation}
% \end{linenomath}
where $U=2\pi r/T$.
The test is run from $t=0$ to $t=2T$.
Along with the initial condition, the state at $t=T$ is shown in the right of Figure \ref{fig:spherical ics}. \\
\\
For a convergence test, the $L^2$ error of $\bm{F}$ is computed for transported solutions at $t=2T$ against the true field, at a range of spatial resolutions.
We took $r=100$ m and $T=200$ s and performed all simulations with $\Delta t=0.05$ s.
To mesh the sphere we use a cubed-sphere grid.
The results for the different schemes are displayed in the right of Figure \ref{fig:convergence transport}, which again shows the improvements of the new schemes of Sections \ref{sec:recovery} and \ref{sec:vorticity} compared with the benchmark scheme of Section \ref{sec:benchmark scheme}.
The results indicate that both schemes are approaching the desired second-order accuracy.

\subsection{Shallow-water test cases} \label{sec:shallow-water}
\begin{figure}[h!]
    \centering
    \includegraphics[width=0.95\textwidth]{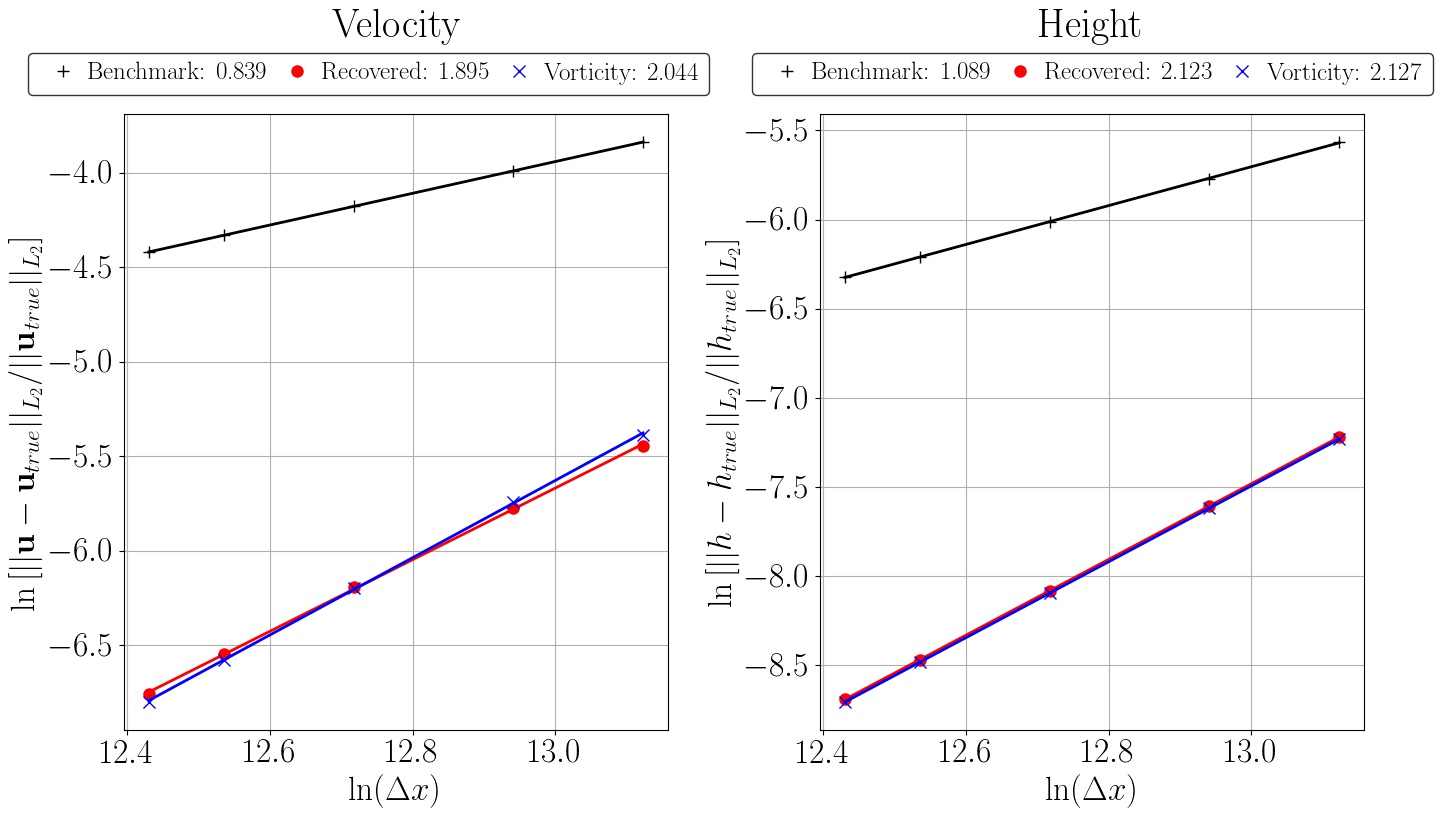}
    \caption{
    Convergence results from the second shallow-water test case test of Williamson et al \cite{williamson1992standard}.
    The normalised $L^2$ error after 5 days is plotted as a function of spatial resolution.
    The test case describes a steady-state zonal geostrophic flow.
    The legends indicate the gradients of lines of best fit through the error measurements, which approximate the order of accuracy of the model.
    The shallow-water simulations differ only in the scheme used to transport the velocity field, comparing the benchmark scheme of Section \ref{sec:benchmark scheme} against the new schemes of Sections \ref{sec:recovery} and \ref{sec:vorticity}.
    (Left) the results for the velocity field $\bm{u}$ and (right) for the height field $h$.
    For both new schemes and for both variables, the model has around second-order accuracy, whereas the benchmark case has only first-order accuracy.
    }
    \label{fig:will2_convergence}
\end{figure}

\begin{figure}[h!]
    \centering
    \includegraphics[width=0.9\textwidth]{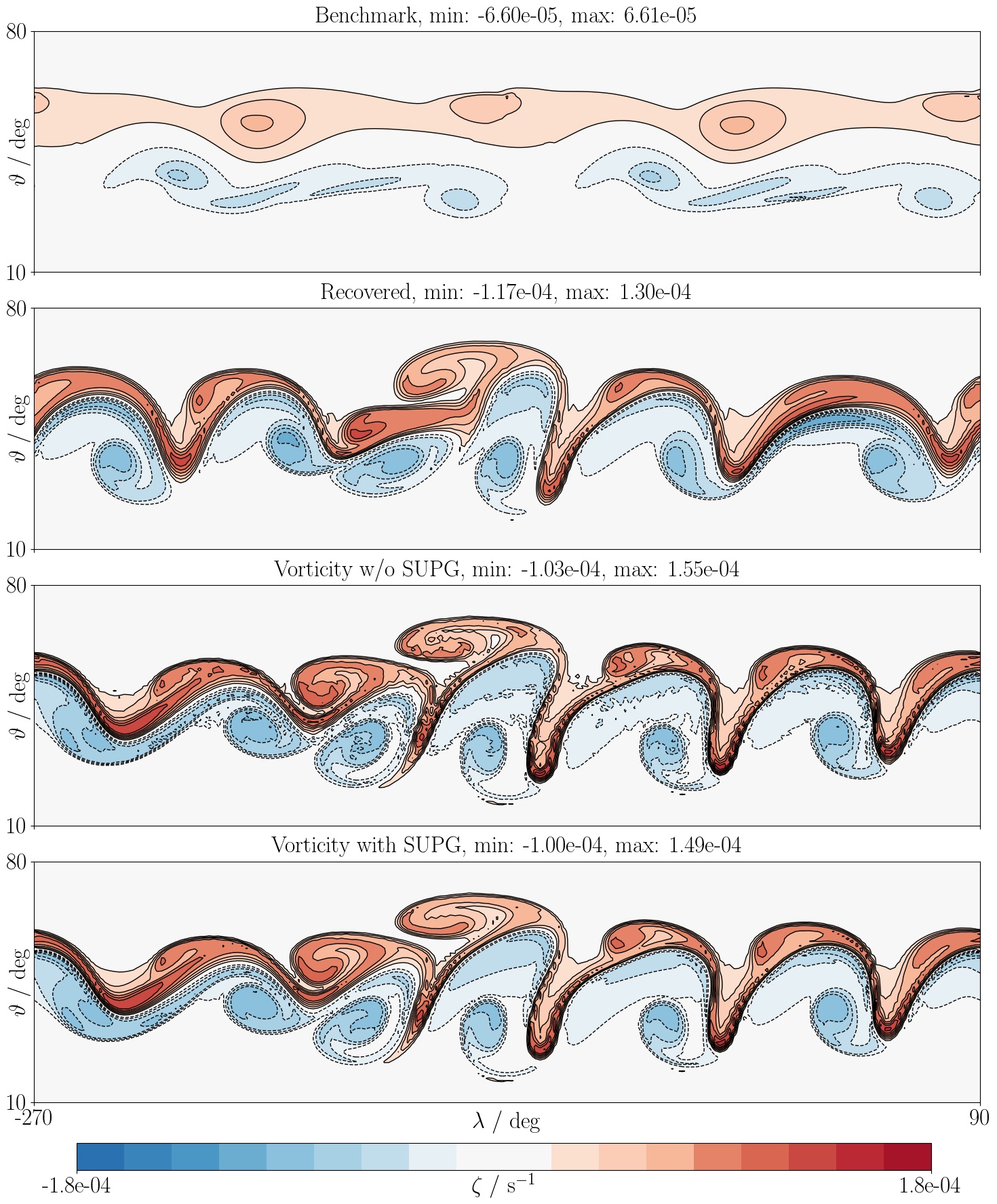}
    \caption{
    The vorticity field after 6 days for the unstable shallow-water jet test case of Galewsky et al \cite{galewsky2004initial}.
    The four plots correspond to simulations which only differ in the velocity transport scheme.
    These simulations were performed on a cubed-sphere mesh with $128\times 128$ cells per panel.
    The simulation using the benchmark scheme is so diffusive that the jet does not clearly form.
    With both the new schemes, the solutions resemble that of \cite{galewsky2004initial}, with the vorticity transport form being particularly close.
    The bottom two panels compare simulations with the vorticity scheme, without and with the SUPG stabilisation.
    The contours are spaced by $2\times 10^{-5}$ s$^{-1}$, and dashed contours indicate negative values.
    The zero contour is omitted.
    }
    \label{fig:galewsky_fields}
\end{figure}

\begin{figure}[h!]
    \centering
    \includegraphics[width=0.95\textwidth]{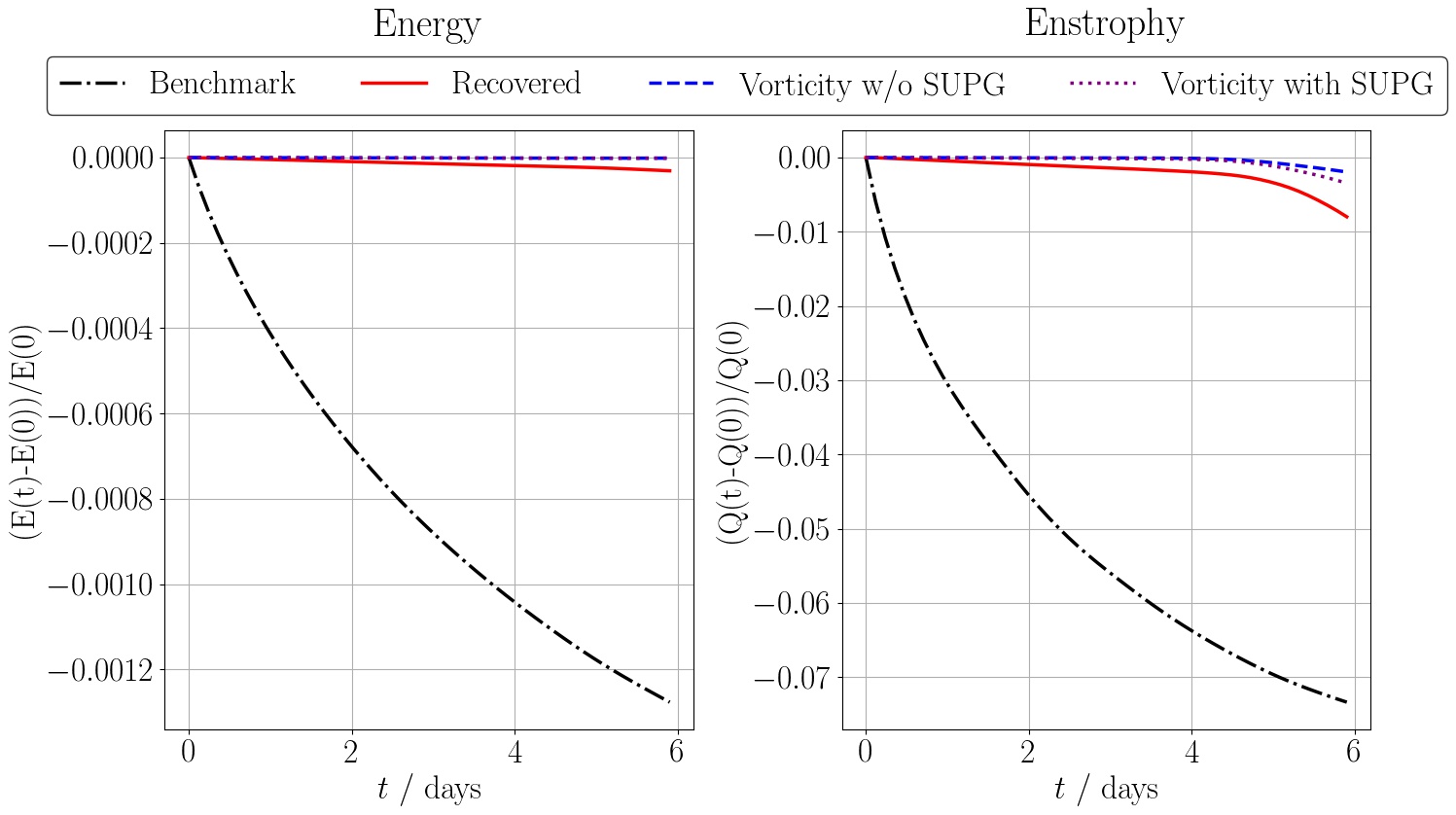}
    \caption{
    Time series of the evolution of the global energy and global enstrophy in the unstable jet test case of \cite{galewsky2004initial}.
    The benchmark scheme shows significant decay of energy and enstrophy compared with the new transport schemes.
    While both vorticity schemes have good conservation of energy, the SUPG stabilisation results in more diffusion of enstrophy.
    }
    \label{fig:galewsky_time_series}
\end{figure}

Now the new vector transport schemes are demonstrated within a compatible finite element discretisation for the shallow-water equations on the sphere \eqref{eqn:shallow-water}.
This discretisation is summarised in Section \ref{sec:shallow-water background}.
The transport of $h$ uses the recovered transport scheme for scalars presented by \cite{bendall2019recovered}, with the time discretisation as the trapezoidal scheme \eqref{eqn:trapezoidal}.
The new schemes for transporting $\bm{u}$ are compared with the benchmark upwind scheme \eqref{eqn:dg upwind}.
For the vorticity scheme in Section \ref{sec:vorticity}, before each transport step the initial vorticity needed for the vorticity evolution equation is updated from $\bm{u}^*$ by solving \eqref{eqn:zeta_def_discr}.
\\
\\
The different velocity transport schemes are demonstrated in this shallow-water model through two standard test cases.
Firstly, the second test from the suite of Williamson et al \cite{williamson1992standard}, which describes a zonal geostrophic flow.
This is a steady-state flow, so the evolved $u$ and $h$ fields can be compared with their initial values to compute errors due to the discretisation.
For full details of the initial conditions, see \cite{williamson1992standard}.
In Figure \ref{fig:will2_convergence}, the errors in $\bm{u}$ and $h$ are plotted after 5 days of simulation, for both the new schemes and the benchmark scheme.
The errors are computed at different spatial resolutions to approximate the order of accuracy of the overall model.
As in Section \ref{sec:sphere test}, the test was performed with a cubed-sphere mesh.
For all simulations we took the same time step of $\Delta t=240$ s.
Figure \ref{fig:will2_convergence}, shows the clear benefits of the two new schemes over the benchmark, improving the order of accuracy of the model from roughly first-order to approximately second-order for both schemes and both variables. \\
 \\
The second test case is the unstable jet of Galewsky et al \cite{galewsky2004initial}.
This test adds a perturbation to an unstable jet in geostrophic balance, which then leads to the jet becoming unbalanced.
Full details of the initial conditions can be found in \cite{galewsky2004initial}.
Figure \ref{fig:galewsky_fields} shows the diagnostic vorticity field after 6 days.
It shows that the benchmark scheme of Section \ref{sec:benchmark scheme} is too diffusive for the fine details of the instability to develop, and that the new schemes are clear improvements on this, with the results resembling those of \cite{galewsky2004initial}.
It also demonstrates the impact of the SUPG stabilisation, by comparing the vorticity scheme with and without this stabilisation (bottom two panels of Figure \ref{fig:galewsky_fields}).
The SUPG stabilisation removes some of the noise seen in the vorticity scheme.
The removal of this noise can also be seen in Figure \ref{fig:galewsky_time_series}, which plots the  evolution over time of the global energy and the global enstrophy for this test case.
While the SUPG stabilisation does not appear to have an effect on the energy, it does result in some degradation of enstrophy compared with the standard vorticity scheme.
Figure \ref{fig:galewsky_time_series} also shows the diffusivity of the benchmark scheme relative to the improved schemes.
These simulations were all performed with a time step of $\Delta t=300$ s, and using a cubed-sphere mesh with $128\times 128$ cells per panel.
\section{Discussion and Summary} \label{sec:summary}
This work has examined two finite element methods for solving the vector transport equation with the lowest-order Raviart-Thomas elements.
This was motivated by increasing the order of accuracy when compared with a standard upwind discretisation.
The first scheme is an extension to the transport schemes of \cite{bendall2019recovered}, and solves the transport equation in advective form and recovers the field in a higher-order function space to transport it there.
The second scheme is a take on the mixed finite element formulation of \cite{bauer2018energy}, applied to the vector transport equation and using a residual based stabilisation concept of \cite{wimmer2020energy}.
This is written in a vorticity form, solving a problem for both the transported vector and its vorticity.
As demonstrated through the test cases in Section \ref{sec:results}, both schemes do have improved accuracy.\\
\\
In the future, we intend to apply these schemes to the lowest-order Raviart-Thomas elements on triangular cells, and to three-dimensional manifolds.
Some preliminary investigations using the test cases described in Section \ref{sec:results} showed that the recovered scheme of Section \ref{sec:recovery} is naturally extended to triangular cells, by using appropriate higher-order finite element spaces for the recovery process.
However, with the lowest-order Raviart-Thomas elements on triangular cells, the vorticity-form scheme of Section \ref{sec:vorticity} suffered from a large amount of noise in the divergence field of $\bm{F}$.
The noise lies in the null-space of the $\bm{\nabla}^\perp\bm{\cdot}$ operator and does not appear in the vorticity evolution equation, and can therefore not be attenuated by the SUPG-based vorticity stabilisation method.
Although it comes at the cost of reduced accuracy, we found that the noise can be controlled effectively by a divergence-based interior penalty term as mentioned in Section \ref{sec:vorticity}.
\section*{Acknowledgements}
The authors would like to thank James Kent, Colin Cotter
and Thomas Melvin for their advice and some very helpful discussions through the evolution of this manuscript.
% \printbibliography

\bibliography{vector_transport}

\end{document}